\documentclass[final,1p,times]{elsarticle}

\usepackage{amsmath,amssymb,amsthm}
\usepackage{booktabs}
\usepackage{algorithm}
\usepackage[noend]{algpseudocode}
\usepackage{placeins}
\usepackage{xcolor}
\usepackage[colorlinks=true,linkcolor=blue,citecolor=blue,urlcolor=blue]{hyperref}
\usepackage[capitalise,nameinlink]{cleveref}
\usepackage{lineno}

\journal{Mathematics and Computers in Simulation}
\newcommand{\Z}{\mathbb{Z}}
\newcommand{\IO}{I_O}
\newcommand{\IT}{I_T}
\newcommand{\Shat}{\hat{S}}
\newcommand{\SNI}{\hat{S}_{\mathrm{NI}}}
\newcommand{\OPT}{\mathrm{OPT}}
\newcommand{\OPTNI}{\mathrm{OPT}_{\mathrm{NI}}}
\newcommand{\Tset}{\mathcal{T}}
\newcommand{\ind}{\mathbf{1}}
\newcommand{\PNI}{\textup{3-PhTP-NI}}
\newcommand{\PTP}{\textup{3-PhTP}}

\theoremstyle{plain}
\newtheorem{theorem}{Theorem}
\newtheorem{lemma}{Lemma}
\newtheorem{proposition}{Proposition}
\newtheorem{corollary}{Corollary}
\theoremstyle{definition}
\newtheorem{definition}{Definition}

\newtheorem{assumption}{Assumption}
\theoremstyle{remark}
\newtheorem{remark}{Remark}

\crefname{assumption}{Assumption}{Assumptions}
\crefname{algorithm}{Algorithm}{Algorithms}
\crefformat{equation}{Eq.~(#2#1#3)}
\Crefformat{equation}{Eq.~(#2#1#3)}
\crefname{figure}{Fig.}{Figs.}

\begin{document}

\begin{frontmatter}

\title{Three-phase time minimization transportation problem with non-idle phases: a polynomial two-witness algorithm}

\author[a,b]{Sailesh Kamalahasan\fnref{eq}}
\author[a,b]{Raghav Paladugu\fnref{eq}}
\author[a]{Archana Khurana\corref{cor1}}
\ead{akhurana@umd.edu}

\cortext[cor1]{Corresponding author.}
\fntext[eq]{These authors contributed equally to this work.}

\affiliation[a]{organization={Department of Mathematics, University of Maryland},
            city={College Park},
            state={MD},
            country={USA}}
\affiliation[b]{organization={Department of Computer Science, University of Maryland},
            city={College Park},
            state={MD},
            country={USA}}

\begin{abstract}
The three-phase time minimization transportation problem partitions the origin--terminal links into three sets used in consecutive phases and minimizes the sum of the phase shipment times. Optimal solutions may leave a phase idle, as in the main example of the original study. We require every phase to carry a positive integer quantity while leaving the split of the total requirement across phases free. Discarding idle solutions from the original algorithm does not, in general, solve this variant. We prove that the optimum is attained by fixing one witness link in each of two phases and deciding the third phase by a single network flow computation, with the optimal Phase-3 time of each witness pair located by binary search. The resulting exact algorithm is strongly polynomial. It extends to any fixed number of phases, whereas, when the number of phases is part of the input, even non-idle feasibility is strongly NP-complete; feasibility alone is fixed-parameter tractable in the number of phases via known rainbow-matching results. We also characterize the minimum time of each phase. On the main example, non-idleness raises the optimum from 38 to 41. On 3,650 random instances with up to $50\times50$ links the algorithm took at most 50.2~s, and non-idleness raised the optimal value in 41\% of the 2,250 instances with equal phase sizes; on 90 sampled instances it was faster than a mixed integer programming solver on 80 with a direct formulation and on 65 with a tighter threshold formulation.
\end{abstract}


\end{frontmatter}

\section{Introduction}\label{sec:intro}

In a time minimization transportation problem (TMTP), a homogeneous commodity is shipped from origins to terminals, all shipments start simultaneously, and the duration of the operation is the largest transportation time among the links that carry a positive quantity. The problem was introduced by Hammer~\cite{Hammer1969}, and Garfinkel and Rao~\cite{Garfinkel1971} studied it as the bottleneck transportation problem. Because the duration depends only on which links are used and not on the amounts carried on them, the TMTP is a natural model for emergency logistics, perishable goods and other settings in which completion time matters more than cost.

Hierarchical and multi-stage versions divide the link set into groups that are used one after another. In the two-level hierarchical TMTP of Sonia and Puri~\cite{Sonia2004}, a quantity is first shipped on Level-I links and the remainder on Level-II links, and the sum of the two shipment times is minimized; Sharma et al.~\cite{Sharma2015} gave a polynomial iterative algorithm for this problem, and Kaur et al.~\cite{Kaur2022} an alternative approach. Capacitated variants with two and three levels have been solved by Xie and Li~\cite{XieLi2022} and by Zhang and Xie~\cite{ZhangXie2024}, and Zhang and Xie~\cite{ZhangXie2025} studied a bilevel model in which the origins are partitioned into two non-empty groups. Two-stage problems with restricted flow~\cite{Kaur2017,Kaur2024rairo} build on the analysis of capacitated transportation problems with bounds on rim conditions by Dahiya and Verma~\cite{Dahiya2007}, and a two-level solid version has also been studied~\cite{KaurJainDahiya2025}.

Jain, Dahiya and Verma~\cite{Jain2021} extended this setting to three phases. In their three-phase time minimization transportation problem (\PTP), the links are partitioned into Phase-1, Phase-2 and Phase-3 links, transportation in each phase starts only after the previous phase is complete, and the sum $T_1+T_2+T_3$ of the three phase shipment times is minimized; their polynomially bounded iterative algorithm generates a sequence of \emph{proficient triads} of phase times. In this journal, Kaur, Dahiya and Sharma~\cite{Kaur2025matcom} studied a two-phase problem in which a specified amount of the commodity is transported in each phase, and a recent preprint by Kaur, Dahiya and Verma~\cite{Kaur2026arxiv} formulates a bi-objective two-phase problem, minimizing the sum of the phase times and the total cost, in which the flow in each phase must be non-zero.

In \PTP{} a phase is allowed to carry nothing, in which case its time is zero. In the main numerical example of~\cite{Jain2021} (their Example~1) the optimal triad is $(06,32,00)$ with value~38, so the optimal schedule transports nothing in Phase-3, and six of the ten proficient triads that the original algorithm records for this example with respect to Phase-3 contain a zero component (\cref{sec:ex1}). An idle phase is harmless when it costs nothing, but not when every phase carries a commitment that does not depend on the amount shipped. Dispatching in waves has been studied, for example, for same-day delivery~\cite{KlappEreraToriello2018}, where whether to dispatch in a given wave is itself a decision. Consider a distributor that dispatches in consecutive waves sharing one fleet: receiving windows and vehicle restrictions determine which routes may run in which wave, and each wave lasts as long as its longest trip. If, in addition, a crew roster or carrier contract requires at least one load per wave, while the total requirement may still be split freely across the waves, no wave may be idle (\cref{sec:dispatch}).

This paper studies the resulting problem, which we call the three-phase time minimization transportation problem with non-idle phases (\PNI): every phase must carry a positive integer quantity, and the split of the total requirement across phases remains a decision, unlike in the restricted-flow model of~\cite{Kaur2025matcom}, which prescribes the amount of each of two phases. \PNI{} is not a routine reformulation of \PTP. The requirement depends on integrality: a real solution may carry positive flow in every phase although no integer solution does (\cref{rem:single,rem:choices}). The auxiliary problems of~\cite{Jain2021} can make a link unusable but cannot require a phase to carry flow; their convergence theory carries over to non-idle solutions (\cref{prop:framework}), but the modified auxiliary problems are no longer standard TMTPs, so the argument of~\cite{Jain2021} for a polynomial bound no longer applies. And discarding the idle triads from the output of the original algorithm does not, in general, give the non-idle optimum: on Example~1 it gives 43, 43 or 41 depending on the phase with respect to which the triads are generated, against the optimum~41 (\cref{sec:ex1}), and in every random sample drawn as in \cref{sec:filterexp} it failed in more than half of the instances, whatever the tie-breaking; the rate depends on how the phases are sampled.

Our method rests on a simple observation. Each phase of a non-idle solution contains a positive link whose time equals the phase time. Fixing such a \emph{witness} link in two phases and reserving one unit on each leaves an ordinary transportation problem on the links whose times do not exceed the witness times, and whether the third phase can then carry positive flow within a given time is decided by one network flow computation (\cref{lem:probe,lem:scc}). Enumerating the witness pairs and locating the optimal Phase-3 time of each pair by a binary search capped by the incumbent value gives an exact algorithm whose number of network flow computations is bounded in terms of the numbers of origins and terminals alone; it is therefore strongly polynomial. The contributions are as follows.
\begin{enumerate}
\item \emph{Model.} We formulate \PNI{} as a single integer transportation problem with three phase-flow constraints and show that non-idleness is equivalent to excluding the triads with a zero component (\cref{sec:model}).
\item \emph{Two-witness reduction and algorithm} (the main result). The optimal value equals the smallest $t_{w_1}+t_{w_2}+\theta_3$ over a Phase-1 link $w_1$, a Phase-2 link $w_2$ and a Phase-3 time value $\theta_3$ for which one auxiliary transportation problem admits a solution with positive Phase-3 flow (\cref{thm:twowitness}). The resulting algorithm is correct (\cref{thm:correct}) and strongly polynomial; it needs at most $N(1+\lceil\log_2\gamma\rceil)$ probes, where $N=\lfloor\frac{mn-1}{2}\rfloor\lceil\frac{mn-1}{2}\rceil$ is the bound stated in~\cite{Jain2021} for the number of its auxiliary problems (a direct count gives at most $N+mn+2$, \cref{sec:theory}), and at most $(mn)^2(1+\lceil\log_2\gamma\rceil)/9$ when the two smallest phases are enumerated (\cref{thm:complexity}). The network flow tools are standard; the reduction is the main contribution of this paper.
\item \emph{Structure of non-idle solutions.} We characterize the minimum attainable time of each phase (\cref{thm:tkmin}), show that it may strictly exceed the smallest time entry of the phase (\cref{prop:strict}), and show by a counterexample that non-idle feasibility cannot be decided from bounds on phase flows alone (\cref{prop:counter}).
\item \emph{Number of phases.} For any fixed number $p$ of phases the reduction gives a strongly polynomial algorithm with at most $\prod_{k<p}|P_k|\,(1+\lceil\log_2\gamma_p\rceil)$ probes (\cref{prop:fixedp}); for $p=2$ it computes the minimum-time point of the bi-objective model of~\cite{Kaur2026arxiv} without costs (\cref{cor:twophase,rem:mincost}). When $p$ is part of the input, deciding whether a non-idle solution exists is strongly NP-complete even with unit availabilities, requirements and times, and minimizing the sum of times is strongly NP-hard even when, in addition, all times lie in $\{1,2\}$ (\cref{thm:hard}); so, unless P${}={}$NP, no algorithm runs in time polynomial in the input size uniformly in $p$, although running times of the form $f(p)\cdot\mathrm{poly}(m,n)$ are not excluded. Feasibility alone is fixed-parameter tractable in $p$ (\cref{prop:fptfeas}); this follows from a simple lemma and known results on rainbow matching~\cite{GuptaRSZ2017} and is not a new technique. We are not aware of an algorithm for the optimization problem with running time $f(p)\cdot\mathrm{poly}(m,n)$, nor of a result that excludes one (\cref{sec:pphases}).
\item \emph{Relation to the original framework.} With three modifications, the convergence theory of~\cite{Jain2021} holds for non-idle solutions (\cref{prop:framework}); the bound $N$ stated there for the number of auxiliary problems can be exceeded, and we give a corrected bound (\cref{sec:theory}).
\item \emph{Computational evidence} (observations, not theorems). We illustrate the method on a synthetic dispatch scenario, a small instance and both examples of~\cite{Jain2021}; for Example~1 the non-idle optimum is $(06,32,03)$ with value~41 (\cref{sec:examples}). A computational study on 3,650 random instances, following the design of~\cite{Kaur2025matcom}, together with further instances for specific questions, shows that the number of probes grows as $(mn)^{1.38}$ when the phases have equal size, that the incumbent cap removes about 77\% of the probes, that non-idleness binds in 41\% of the instances with equal phase sizes, and that the algorithm is faster than a mixed integer programming solver at the median, by a factor of 2.8 with a direct formulation and 1.5 with a tighter threshold formulation (\cref{sec:experiments}). These findings depend on the instance distributions and implementations used.
\end{enumerate}

The closest prior work is the preprint~\cite{Kaur2026arxiv}. Its feasible sets require each of two phases to carry a non-zero integer flow, and by its Theorem~2 the last non-dominated point recorded by its algorithm attains the minimum sum of the two phase times over such solutions. The present paper treats three phases, and more generally any fixed number of phases, with the single objective of the total time; builds the requirement into the algorithm through witnesses covered by the correctness proof (\cref{thm:correct}); and establishes the complexity of the problem in the number of phases. As a secondary observation, we show in \cref{sec:othermodels} that the re-solution rule by which~\cite{Kaur2026arxiv} enforces the requirement (its Remark~12) can discard every non-idle solution of a sub-problem, after which its stopping rule certifies a sum of times that is not minimal.

\Cref{sec:model} states the model, \cref{sec:algorithm} presents the algorithm, and \cref{sec:theory} establishes its correctness and complexity and the results on phase minima, feasibility and the number of phases. \Cref{sec:examples,sec:experiments} report numerical illustrations and computational experiments, \cref{sec:discussion} discusses the relation to the framework of~\cite{Jain2021} and to related models, and \cref{sec:conclusions} concludes. The Supplementary Material contains additional tables, the full specification of the sampling study, the HiGHS reproduction, two further instances for \cref{sec:othermodels}, computational checks of the proofs, and the scripts and data.

\section{Problem formulation}\label{sec:model}

\subsection{Data and notation}

Let $\IO=\{1,\dots,m\}$ and $\IT=\{1,\dots,n\}$ be the index sets of origins and terminals. Origin $i$ has availability $a_i$ and terminal $j$ has requirement $b_j$, where $a_i,b_j$ are positive integers and the problem is balanced, $\sum_{i\in\IO}a_i=\sum_{j\in\IT}b_j=:Q$. The time of transporting the commodity along link $(i,j)$ is $t_{ij}$, independent of the quantity shipped. The link set $\IO\times\IT$ is partitioned into three non-empty disjoint sets $P_1$, $P_2$ and $P_3$ of Phase-1, Phase-2 and Phase-3 links.

\begin{assumption}\label{ass:positive}
All transportation times are positive: $t_{ij}>0$ for every $(i,j)\in\IO\times\IT$.
\end{assumption}

For an allocation $X=(x_{ij})$ we define, for $k=1,2,3$, the Phase-$k$ flow and the Phase-$k$ shipment time
\begin{equation}\label{eq:FT}
\begin{aligned}
F_k(X)&=\sum_{(i,j)\in P_k}x_{ij},\\
T_k(X)&=\begin{cases}\max\{t_{ij}:(i,j)\in P_k,\ x_{ij}>0\}, & \text{if } F_k(X)>0,\\ 0, & \text{if } F_k(X)=0.\end{cases}
\end{aligned}
\end{equation}
The triple $(T_1(X),T_2(X),T_3(X))$ is the \emph{triad} of $X$. We write $\tau_k=\min\{t_{ij}:(i,j)\in P_k\}$ for the smallest Phase-$k$ time, $\Tset_k$ for the set of distinct Phase-$k$ times, and $\gamma=|\Tset_3|$. For a link $w=(i,j)$ we write $t_w=t_{ij}$ and $\ind_w$ for the $m\times n$ matrix with a one in position $w$ and zeros elsewhere.

\subsection{The model}

Let
\begin{equation}\label{eq:S}
\Shat=\Big\{X\in\Z^{m\times n}:\ \sum_{j\in\IT}x_{ij}=a_i\ (i\in\IO),\ \sum_{i\in\IO}x_{ij}=b_j\ (j\in\IT),\ x_{ij}\ge 0\Big\}
\end{equation}
be the set of integer feasible solutions of the balanced transportation problem. \PTP{} is the problem $\min\{T_1(X)+T_2(X)+T_3(X):X\in\Shat\}$ with optimal value $\OPT$~\cite{Jain2021}. The set of \emph{non-idle} solutions is
\begin{equation}\label{eq:SNI}
\SNI=\{X\in\Shat:\ F_k(X)\ge 1,\ k=1,2,3\},
\end{equation}
and the problem studied in this paper is
\begin{equation}\label{eq:NI}
\PNI:\qquad \OPTNI=\min_{X\in\SNI}\ T_1(X)+T_2(X)+T_3(X),
\end{equation}
with $\OPTNI=+\infty$ if $\SNI=\emptyset$.

\begin{remark}[Single-matrix formulation]\label{rem:single}
In~\cite{Jain2021} the problem is stated with three real allocation matrices $Y$, $Z$ and $W$ for the three phases. The availabilities and requirements left after Phase-1 and Phase-2 define the feasible sets of the later phases, and the matrices are supported on $P_1$, $P_2$ and $P_3$, respectively (the condition on $W$ printed in~\cite{Jain2021}, $w_{ij}=0$ for $(i,j)\in P_2\cup P_3$, must be read with $P_1\cup P_2$). Since the phase link sets are disjoint, $X=Y+Z+W$ maps these feasible triples one-to-one onto the real transportation polytope, and maps the integer triples one-to-one onto $\Shat$; the inverse restricts $X$ to $P_1$, $P_2$ and $P_3$. The map preserves the phase allocations, the phase flows $F_k$, the phase times $T_k$ and hence the objective, so the residual availabilities and requirements need not be modeled separately.

The optimal value of \PTP{} is the same over the real polytope and over $\Shat$. Indeed, the solutions supported in the support of a real solution form a non-empty face of the polytope, and a vertex of this face is an integral vertex of the polytope that uses no new link. The auxiliary problems of~\cite{Jain2021} are posed over a single matrix as well, and the same argument shows that their optimal values are attained in $\Shat$.

Non-idleness behaves differently: in the instance of \cref{prop:counter} a real solution has positive flow in every phase, although $\SNI=\emptyset$ (\cref{rem:choices}).
\end{remark}

\begin{remark}[Modeling choices]\label{rem:choices}
Two choices in~\cref{eq:SNI} are deliberate. First, the model is posed over integers. Over the reals, non-idleness would be the strict inequality $F_k(X)>0$, and the feasible set would change: in the instance of the proof of \cref{prop:counter} every $X$ with $x_{11}=s\in(0,1)$ is a real non-idle solution, whereas $\SNI=\emptyset$. (Since $T_1+T_2+T_3$ depends only on the support of $X$, it takes finitely many values, so a minimum over real non-idle solutions would still be attained whenever such solutions exist.) Over the integers, $F_k(X)>0$ is equivalent to $F_k(X)\ge1$. Second, the constraint is a positivity requirement and not a prescribed quantity $F_k(X)\ge L_k$ or $F_k(X)=L_k$. Prescribing the quantity of each phase, as in the restricted-flow model of~\cite{Kaur2025matcom}, fixes the split of the requirement across phases; in \PNI{} the split is part of the decision.
\end{remark}

The next lemma shows that, under \cref{ass:positive}, non-idleness is a condition on the triad alone.

\begin{lemma}\label{lem:idle}
Under \cref{ass:positive}, for every $X\in\Shat$ and $k\in\{1,2,3\}$,
\[
F_k(X)\ge1\iff T_k(X)>0\iff T_k(X)\ge\tau_k .
\]
\end{lemma}

\begin{proof}
If $F_k(X)\ge1$, some link $(i,j)\in P_k$ has $x_{ij}>0$, so $T_k(X)\ge t_{ij}\ge\tau_k>0$. If $F_k(X)=0$, then $T_k(X)=0$ by~\cref{eq:FT}. Hence $F_k(X)\ge1$ implies $T_k(X)\ge\tau_k$, which implies $T_k(X)>0$ because $\tau_k>0$; and $T_k(X)>0$ implies $F_k(X)\neq0$, that is, $F_k(X)\ge1$, because $F_k(X)$ is a non-negative integer. So the three conditions are equivalent.
\end{proof}

Consequently, $\SNI$ consists exactly of the solutions of \PTP{} whose triad has no zero component, and \PNI{} can be viewed as \PTP{} with the triads containing a zero removed. The following facts are immediate.

\begin{proposition}\label{prop:relation}
$\OPT\le\OPTNI$. If some optimal solution of \PTP{} belongs to $\SNI$, then $\OPTNI=\OPT$.
\end{proposition}

\begin{proof}
The first claim follows from $\SNI\subseteq\Shat$. If an optimal solution $X^*$ of \PTP{} is non-idle, it is feasible for~\cref{eq:NI}, so $\OPTNI\le T_1(X^*)+T_2(X^*)+T_3(X^*)=\OPT$. Together with the first claim, $\OPTNI=\OPT$.
\end{proof}

\Cref{prop:relation} does not apply to Example~1 of~\cite{Jain2021}, whose optimal triad $(06,32,00)$ is idle; we show in \cref{sec:ex1} that $\OPTNI=41>38=\OPT$ for that instance.

\section{The two-witness algorithm}\label{sec:algorithm}

\subsection{Witnesses and the auxiliary problem}

The algorithm fixes one link in Phase-1 and one link in Phase-2 that are required to carry at least one unit; these links are called \emph{witnesses}. Reserving the two units reduces the availabilities and requirements, and the remaining quantity is shipped on links whose times do not exceed the witness times.

\begin{definition}\label{def:aux}
Let $w_1\in P_1$, $w_2\in P_2$ and $\theta\in\Tset_3$. For $i\in\IO$ and $j\in\IT$ let $\rho_i$ and $\sigma_j$ be the number of witnesses in row $i$ and column $j$, respectively, and define the reduced margins $a'_i=a_i-\rho_i$ and $b'_j=b_j-\sigma_j$. The \emph{admissible link set} is
\begin{equation}\label{eq:L}
\begin{aligned}
L(w_1,w_2,\theta)={}&\{(i,j)\in P_1:t_{ij}\le t_{w_1}\}\cup\{(i,j)\in P_2:t_{ij}\le t_{w_2}\}\\
&\cup\{(i,j)\in P_3:t_{ij}\le\theta\}.
\end{aligned}
\end{equation}
The \emph{auxiliary problem} $R(w_1,w_2,\theta)$ is the transportation problem with margins $a',b'$ in which $x_{ij}=0$ for every $(i,j)\notin L(w_1,w_2,\theta)$. We say that \emph{$\Pi(w_1,w_2,\theta)$ holds} if all reduced margins are non-negative, $R(w_1,w_2,\theta)$ has an integer feasible solution, and some integer feasible solution $X'$ of $R(w_1,w_2,\theta)$ satisfies $F_3(X')\ge1$.
\end{definition}

\subsection{Evaluating \texorpdfstring{$\Pi$}{Pi} by one minimum cost flow}\label{sec:probe}

Whether $\Pi(w_1,w_2,\theta)$ holds can be decided by solving one cost-minimizing transportation problem. Let $Q'=Q-2$ be the total reduced availability, let $M$ be any number with $M>Q'$, and define the costs
\begin{equation}\label{eq:probecost}
c_{ij}=\begin{cases}
M, & (i,j)\notin L(w_1,w_2,\theta),\\
-1, & (i,j)\in L(w_1,w_2,\theta)\cap P_3,\\
0, & (i,j)\in L(w_1,w_2,\theta)\setminus P_3 .
\end{cases}
\end{equation}
We call the problem $\min\{\sum_{i,j}c_{ij}x_{ij}\}$ over the transportation polytope with margins $a',b'$ the \emph{probe} for $(w_1,w_2,\theta)$.

\begin{lemma}\label{lem:probe}
Assume $a'\ge0$ and $b'\ge0$, and let $X'$ be an optimal basic solution of the probe. Then $X'$ is integral, and $\Pi(w_1,w_2,\theta)$ holds if and only if $x'_{ij}=0$ for every $(i,j)\notin L(w_1,w_2,\theta)$ and $\sum_{i,j}c_{ij}x'_{ij}\le-1$. In that case $X'$ is feasible for $R(w_1,w_2,\theta)$ and $F_3(X')$ equals the largest Phase-3 flow over all feasible solutions of $R(w_1,w_2,\theta)$.
\end{lemma}

\begin{proof}
The transportation polytope with margins $a',b'\ge0$ is non-empty, because $\sum_ia'_i=\sum_jb'_j=Q'$, and bounded, so the probe has an optimal basic solution, and every basic solution is integral by the total unimodularity of the transportation constraint matrix and the integrality of $a',b'$. Let $Y$ be an integral point of the polytope; then $\sum_{i,j}y_{ij}=Q'$. If $Y$ carries at least one unit on a link outside $L$, its cost is at least $M-(Q'-1)>0$. If $Y$ vanishes outside $L$, it is an integer feasible solution of $R(w_1,w_2,\theta)$ and its cost is $-F_3(Y)\le0$; more generally, every real feasible solution $Z$ of $R(w_1,w_2,\theta)$ is a point of the polytope with cost $-F_3(Z)$.

If $\Pi(w_1,w_2,\theta)$ holds, let $X''$ be an integer feasible solution of $R(w_1,w_2,\theta)$ with $F_3(X'')\ge1$. By optimality, $\sum c_{ij}x'_{ij}\le-F_3(X'')\le-1<0$; since $X'$ is integral, it cannot carry flow outside $L$, so it is feasible for $R(w_1,w_2,\theta)$, and its cost is at most $-1$. Conversely, if $X'$ vanishes outside $L$ and has cost at most $-1$, then $X'$ is an integer feasible solution of $R(w_1,w_2,\theta)$ with $F_3(X')=-\sum c_{ij}x'_{ij}\ge1$, so $\Pi(w_1,w_2,\theta)$ holds. In that case, for every feasible solution $Z$ of $R(w_1,w_2,\theta)$, $-F_3(X')=\sum c_{ij}x'_{ij}\le\sum c_{ij}z_{ij}=-F_3(Z)$, so $F_3(X')$ is the largest Phase-3 flow over all feasible solutions of $R(w_1,w_2,\theta)$.
\end{proof}

In hand computations the probe can be solved by the northwest corner rule followed by the MODI (u--v) method, with $M$ kept symbolic; in implementations any minimum cost flow algorithm may be used. The links outside $L$ can equally be deleted from the network, in which case infeasibility is reported by the flow algorithm directly.

The probe can also be evaluated without costs, by one maximum flow computation and one search for strongly connected components. For a feasible solution $X^0$ of $R(w_1,w_2,\theta)$, let $D(X^0)$ be the directed graph on $\IO\cup\IT$ with an arc $i\to j$ for every $(i,j)\in L(w_1,w_2,\theta)$ and an arc $j\to i$ for every $(i,j)$ with $x^0_{ij}>0$.

\begin{lemma}\label{lem:scc}
Assume $a'\ge0$, $b'\ge0$, and let $X^0$ be an integer feasible solution of $R(w_1,w_2,\theta)$. Then $\Pi(w_1,w_2,\theta)$ holds if and only if $F_3(X^0)\ge1$ or some link $(i,j)\in L(w_1,w_2,\theta)\cap P_3$ has $i$ and $j$ in the same strongly connected component of $D(X^0)$.
\end{lemma}

\begin{proof}
($\Leftarrow$) If $F_3(X^0)\ge1$ there is nothing to prove. Otherwise let $(i,j)\in L\cap P_3$ with $i,j$ in one strongly connected component; then $x^0_{ij}=0$, so $D(X^0)$ has no arc $j\to i$. Take a shortest directed path from $j$ to $i$ in $D(X^0)$. Its arcs alternate between arcs $j'\to i'$ with $x^0_{i'j'}\ge1$ and arcs $i'\to j'$ with $(i',j')\in L$, and it has at least three arcs. Closing it with the arc $i\to j$ gives a simple directed cycle with at least four arcs, so no link occurs on it both as an arc $i'\to j'$ and as an arc $j'\to i'$. Adding one unit on the links of the arcs $i'\to j'$ and subtracting one unit on the links of the arcs $j'\to i'$ preserves every row and column sum, because each node of the cycle has one arc of each kind, and preserves non-negativity, because $X^0$ is integral. The support stays in $L$, and the resulting $X$ has $x_{ij}=1$. It is feasible for $R(w_1,w_2,\theta)$ with $F_3(X)\ge1$.

($\Rightarrow$) Let $X$ be an integer feasible solution of $R(w_1,w_2,\theta)$ with $x_{ij}\ge1$ for some $(i,j)\in P_3$; then $(i,j)\in L$. If $x^0_{ij}\ge1$, then $F_3(X^0)\ge1$. Otherwise suppose, for a contradiction, that $i$ is not reachable from $j$, and let $S$ be the set of nodes reachable from $j$ in $D(X^0)$, so that $j\in S$ and $i\notin S$. The matrix $\Delta=X-X^0$ has zero row and column sums, hence
\[
0=\sum_{i'\in S\cap\IO}\ \sum_{j'\in\IT}\Delta_{i'j'}-\sum_{j'\in S\cap\IT}\ \sum_{i'\in\IO}\Delta_{i'j'}
=\sum_{i'\in S,\ j'\notin S}\Delta_{i'j'}-\sum_{i'\notin S,\ j'\in S}\Delta_{i'j'} .
\]
If $i'\in S$, $j'\notin S$ and $\Delta_{i'j'}>0$, then $x_{i'j'}>0$, so $(i',j')\in L$ and $j'$ is reachable from $j$; hence $\Delta_{i'j'}\le0$ in the first sum. If $i'\notin S$, $j'\in S$ and $\Delta_{i'j'}<0$, then $x^0_{i'j'}>0$, so $D(X^0)$ has the arc $j'\to i'$ and $i'\in S$; hence $\Delta_{i'j'}\ge0$ in the second sum. The identity therefore forces every term of the second sum to vanish, but the term with $i'=i$, $j'=j$ is $\Delta_{ij}=x_{ij}\ge1$. So $j$ reaches $i$, and with the arc $i\to j$ the two nodes lie in one strongly connected component.
\end{proof}

The first alternative in \cref{lem:scc} implies the second and can therefore be omitted: if $x^0_{ij}>0$ for some $(i,j)\in P_3$, then $(i,j)\in L$, and the arcs $i\to j$ and $j\to i$ put $i$ and $j$ in one strongly connected component.

A feasible $X^0$ is an integral maximum flow on the bipartite network with $m+n+2$ nodes and at most $mn+m+n$ arcs, and the strongly connected components are found in linear time~\cite{Tarjan1972}. The computations of \cref{sec:experiments} use this evaluation; they were cross-checked against the minimum cost flow evaluation of \cref{lem:probe} on 20,000 random probes without disagreement.

\subsection{The algorithm}

The algorithm enumerates witness pairs in order of increasing times and, for each pair, finds the smallest $\theta\in\Tset_3$ for which $\Pi$ holds, whenever that value can improve the incumbent value $z^{\mathrm{inc}}$. By \cref{lem:monotone} below, $\Pi$ is monotone in $\theta$. Only values $\theta<z^{\mathrm{inc}}-t_{w_1}-t_{w_2}$ can improve the incumbent, so the algorithm first tests the largest such value, $\theta^h$; if $\Pi$ fails there, the pair is discarded after a single probe, and otherwise the smallest admissible value is located by binary search on $\theta^1,\dots,\theta^h$. Two pruning rules use the incumbent: no pair can improve on $z^{\mathrm{inc}}$ unless $t_{w_1}+t_{w_2}+\tau_3<z^{\mathrm{inc}}$.

\begin{algorithm}[htbp]
\caption{Two-witness algorithm for \PNI}\label{alg:tw}
\begin{algorithmic}[1]
\Require Availabilities $a$, requirements $b$, times $t$, partition $P_1,P_2,P_3$ satisfying \cref{ass:positive}.
\Ensure $\OPTNI$ and an optimal solution $X^*\in\SNI$, or the statement $\SNI=\emptyset$.
\State Relabel the phases so that $|P_3|\ge\max\{|P_1|,|P_2|\}$ \Comment{\cref{rem:symmetry}}
\State Sort $P_1=(w_1^{(1)},\dots,w_1^{(|P_1|)})$ and $P_2=(w_2^{(1)},\dots,w_2^{(|P_2|)})$ by non-decreasing time; let $\Tset_3=\{\theta^1<\dots<\theta^\gamma\}$, $\tau_2=t_{w_2^{(1)}}$, $\tau_3=\theta^1$.
\State $z^{\mathrm{inc}}\gets+\infty$; $X^*\gets$ undefined
\For{$p=1,\dots,|P_1|$}
  \State $w_1\gets w_1^{(p)}$
  \If{$t_{w_1}+\tau_2+\tau_3\ge z^{\mathrm{inc}}$} \textbf{stop the outer loop} \Comment{outer pruning}
  \EndIf
  \For{$q=1,\dots,|P_2|$}
    \State $w_2\gets w_2^{(q)}$
    \If{$t_{w_1}+t_{w_2}+\tau_3\ge z^{\mathrm{inc}}$} \textbf{break} \Comment{inner pruning}
    \EndIf
    \State Compute $a',b'$; \textbf{if} some $a'_i<0$ or $b'_j<0$ \textbf{then continue}
    \State $h\gets\max\{i:\theta^i<z^{\mathrm{inc}}-t_{w_1}-t_{w_2}\}$ \Comment{incumbent cap; $h\ge1$}
    \If{$\Pi(w_1,w_2,\theta^h)$ does not hold} \textbf{continue} \Comment{one probe}
    \EndIf
    \State $lo\gets1$; $hi\gets h$ \Comment{invariant: $\Pi(w_1,w_2,\theta^{hi})$ holds}
    \While{$lo<hi$}
      \State $mid\gets\lfloor(lo+hi)/2\rfloor$
      \If{$\Pi(w_1,w_2,\theta^{mid})$ holds} $hi\gets mid$ \Else\ $lo\gets mid+1$
      \EndIf
    \EndWhile
    \State $z^{\mathrm{inc}}\gets t_{w_1}+t_{w_2}+\theta^{hi}$ \Comment{a strict improvement}
    \State $X^*\gets X'+\ind_{w_1}+\ind_{w_2}$, where $X'$ is the solution of $R(w_1,w_2,\theta^{hi})$ with $F_3(X')\ge1$ found by the last probe in which $\Pi(w_1,w_2,\theta^{hi})$ held (line~12 or~16)
  \EndFor
\EndFor
\If{$z^{\mathrm{inc}}=+\infty$} \Return ``$\SNI=\emptyset$''
\Else\ \Return $\OPTNI=z^{\mathrm{inc}}$ and $X^*$
\EndIf
\end{algorithmic}
\end{algorithm}

When $z^{\mathrm{inc}}=+\infty$ the cap gives $h=\gamma$. The bound $h\ge1$ holds because the inner pruning test failed, so $\theta^1=\tau_3<z^{\mathrm{inc}}-t_{w_1}-t_{w_2}$.
The solution $X'$ used in line~19 is available without a further probe: when probes are evaluated by \cref{lem:probe}, it is the optimal basic solution of the probe; when they are evaluated by \cref{lem:scc}, it is $X^0$ if $F_3(X^0)\ge1$, and otherwise $X^0$ modified along the cycle constructed in the proof of \cref{lem:scc}.

\begin{remark}[Choice of the enumerated phases]\label{rem:symmetry}
The objective in~\cref{eq:NI} and the constraints in~\cref{eq:SNI} are symmetric in the three phase labels; the order in which phases are executed does not affect them. \Cref{alg:tw} may therefore be applied with any two phases as the enumerated ones and the remaining phase in the role of Phase-3. Enumerating the two phases with the fewest links, and letting the phase with the most links be decided by the probe, minimizes the number of witness pairs.
\end{remark}

\FloatBarrier
\section{Theoretical analysis}\label{sec:theory}

\subsection{Correctness}

\begin{theorem}[Two-witness reduction]\label{thm:twowitness}
Under \cref{ass:positive},
\begin{equation}\label{eq:twowitness}
\OPTNI=\min\big\{t_{w_1}+t_{w_2}+\theta\ :\ w_1\in P_1,\ w_2\in P_2,\ \theta\in\Tset_3,\ \Pi(w_1,w_2,\theta)\text{ holds}\big\},
\end{equation}
where the minimum of an empty set is $+\infty$. In particular, $\SNI\neq\emptyset$ if and only if $\Pi(w_1,w_2,\max\Tset_3)$ holds for some $w_1\in P_1$ and $w_2\in P_2$.
\end{theorem}

\begin{proof}
($\le$) Let $\Pi(w_1,w_2,\theta)$ hold, with $X'$ an integer feasible solution of $R(w_1,w_2,\theta)$ satisfying $F_3(X')\ge1$, and set $X=X'+\ind_{w_1}+\ind_{w_2}$. Adding the two units restores the original margins, so $X\in\Shat$. We have $F_1(X)\ge x_{w_1}\ge1$, $F_2(X)\ge x_{w_2}\ge1$ and $F_3(X)=F_3(X')\ge1$, so $X\in\SNI$. Every positive entry of $X$ lies in $L(w_1,w_2,\theta)$, because $w_1$ and $w_2$ belong to it. Hence $T_1(X)\le t_{w_1}$, and since $x_{w_1}\ge1$, in fact $T_1(X)=t_{w_1}$; likewise $T_2(X)=t_{w_2}$, and $T_3(X)\le\theta$. Therefore $\OPTNI\le t_{w_1}+t_{w_2}+\theta$.

($\ge$) If $\SNI=\emptyset$ there is nothing to prove. Otherwise, since $\Shat$ is finite, \cref{eq:NI} has an optimal solution $X^*$; let $(T_1,T_2,T_3)$ be its triad. As $F_k(X^*)\ge1$, each phase carries a positive entry. Choose $w_1\in P_1$ with $x^*_{w_1}>0$ and $t_{w_1}=T_1$, choose $w_2\in P_2$ with $x^*_{w_2}>0$ and $t_{w_2}=T_2$, and put $\theta=T_3\in\Tset_3$. Then $X'=X^*-\ind_{w_1}-\ind_{w_2}$ is a non-negative integer matrix with margins $a',b'$, so the reduced margins are non-negative. Every positive entry of $X^*$ in $P_k$ has time at most $T_k$, so $X'$ is supported on $L(w_1,w_2,\theta)$ and is feasible for $R(w_1,w_2,\theta)$. Finally $F_3(X')=F_3(X^*)\ge1$, so $\Pi(w_1,w_2,\theta)$ holds and the right-hand side of~\cref{eq:twowitness} is at most $T_1+T_2+T_3=\OPTNI$. The two inequalities together prove~\cref{eq:twowitness}.

For the last statement, $\SNI\neq\emptyset$ if and only if $\OPTNI<+\infty$, since $\OPTNI$ is a minimum over the finite set $\SNI$ and $\OPTNI=+\infty$ if $\SNI=\emptyset$. By~\cref{eq:twowitness}, this holds if and only if $\Pi(w_1,w_2,\theta)$ holds for some $w_1\in P_1$, $w_2\in P_2$ and $\theta\in\Tset_3$, and by \cref{lem:monotone} this is the case if and only if $\Pi(w_1,w_2,\max\Tset_3)$ holds for some $w_1\in P_1$ and $w_2\in P_2$.
\end{proof}

\begin{lemma}[Monotonicity]\label{lem:monotone}
If $\Pi(w_1,w_2,\theta)$ holds and $\theta\le\theta'\in\Tset_3$, then $\Pi(w_1,w_2,\theta')$ holds.
\end{lemma}

\begin{proof}
The reduced margins do not depend on $\theta$, and $L(w_1,w_2,\theta)\subseteq L(w_1,w_2,\theta')$ by~\cref{eq:L}. A feasible solution of $R(w_1,w_2,\theta)$ with positive Phase-3 flow is therefore also feasible for $R(w_1,w_2,\theta')$ with the same Phase-3 flow.
\end{proof}

\begin{theorem}[Correctness of \cref{alg:tw}]\label{thm:correct}
Under \cref{ass:positive}, \cref{alg:tw} terminates. If $\SNI=\emptyset$ it reports this, and otherwise it returns $\OPTNI$ and a solution $X^*\in\SNI$ whose objective value is $\OPTNI$.
\end{theorem}

\begin{proof}
Termination is clear, since both loops are finite and each binary search shrinks the interval $[lo,hi]$. For a fixed pair $(w_1,w_2)$ with non-negative reduced margins, let $\theta^\star(w_1,w_2)$ denote the smallest $\theta\in\Tset_3$ for which $\Pi$ holds, if any. By \cref{lem:monotone}, $\Pi(w_1,w_2,\theta)$ holds exactly for $\theta\ge\theta^\star$. A pair with a negative reduced margin admits no $X\in\Shat$ with $x_{w_1}\ge1$ and $x_{w_2}\ge1$, so skipping it loses nothing. By \cref{thm:twowitness}, $\OPTNI$ is the minimum of $t_{w_1}+t_{w_2}+\theta^\star(w_1,w_2)$ over all pairs for which $\theta^\star$ exists.

It remains to show that neither the pruning rules nor the cap discard a pair, or a value of $\theta$, that could improve the incumbent. Every candidate value of a pair $(w_1,w_2)$ is at least $t_{w_1}+t_{w_2}+\tau_3$. In the inner loop the pairs share $w_1$ and have non-decreasing $t_{w_2}$, so once $t_{w_1}+t_{w_2}+\tau_3\ge z^{\mathrm{inc}}$ holds it holds for all later $w_2$, and breaking is safe. In the outer loop every pair with the current or a later $w_1$ has value at least $t_{w_1}+\tau_2+\tau_3$, so stopping is safe when this bound reaches $z^{\mathrm{inc}}$.

Consider a pair that reaches the cap. A value $\theta$ yields a candidate below $z^{\mathrm{inc}}$ only if $\theta<z^{\mathrm{inc}}-t_{w_1}-t_{w_2}$, that is, only if $\theta\le\theta^h$. If $\Pi(w_1,w_2,\theta^h)$ fails, \cref{lem:monotone} shows that $\Pi(w_1,w_2,\theta)$ fails for every $\theta\le\theta^h$, so the pair has no candidate below $z^{\mathrm{inc}}$ and skipping it is safe. Otherwise $\theta^\star(w_1,w_2)\le\theta^h$, the binary search maintains that $\Pi$ holds at $\theta^{hi}$ and fails below $\theta^{lo}$, and it ends with $\theta^{hi}=\theta^\star(w_1,w_2)$; the new value $t_{w_1}+t_{w_2}+\theta^\star\le t_{w_1}+t_{w_2}+\theta^h<z^{\mathrm{inc}}$ is the pair's best candidate. Hence, at termination, $z^{\mathrm{inc}}$ equals the minimum in \cref{eq:twowitness} over all pairs, which is $\OPTNI$ by \cref{thm:twowitness}. In particular, $z^{\mathrm{inc}}=+\infty$ at termination if and only if $\SNI=\emptyset$, and in that case \cref{alg:tw} reports $\SNI=\emptyset$. When the incumbent is updated, the recorded solution $X^*=X'+\ind_{w_1}+\ind_{w_2}$ lies in $\SNI$ and, by the first part of the proof of \cref{thm:twowitness}, has objective value at most $t_{w_1}+t_{w_2}+\theta^\star=z^{\mathrm{inc}}$. Since its value is also at least $\OPTNI$, the solution returned at termination is optimal.
\end{proof}

\subsection{Complexity}

\begin{theorem}[Complexity]\label{thm:complexity}
Let $e$ be the number of witness pairs evaluated by \cref{alg:tw} (pairs on which at least one probe is solved) and $I$ the number of incumbent updates. \Cref{alg:tw} solves at most $e+I\lceil\log_2\gamma\rceil$ probes. Moreover $I\le e\le|P_1||P_2|$, so at most
\begin{equation}\label{eq:N}
|P_1|\,|P_2|\,(1+\lceil\log_2\gamma\rceil)\ \le\ N\,(1+\lceil\log_2\gamma\rceil),\qquad N=\Big\lfloor\frac{mn-1}{2}\Big\rfloor\Big\lceil\frac{mn-1}{2}\Big\rceil,
\end{equation}
probes are solved, and $|P_1||P_2|\le(mn)^2/9$ after the relabelling of line~1. If all times are integers, $I\le\sum_{k=1}^3\big(\max_{(i,j)\in P_k}t_{ij}-\tau_k\big)+1$. Each probe takes $O\big((m+n)\log(m+n)\,(mn+(m+n)\log(m+n))\big)$ time when evaluated by \cref{lem:probe} with Orlin's algorithm~\cite{Orlin1993}, and $O\big((m+n)\,mn\big)$ time when evaluated by \cref{lem:scc} with an $O(VE)$ maximum flow algorithm~\cite{Orlin2013}. \Cref{alg:tw} therefore runs in strongly polynomial time, $O\big((mn)^3(m+n)\log(mn)\big)$ with the latter evaluation.
\end{theorem}

\begin{proof}
An evaluated pair costs one probe at $\theta^h$. If $\Pi(w_1,w_2,\theta^h)$ fails, nothing more is spent on the pair. If it holds, the binary search on an interval of $h\le\gamma$ values costs at most $\lceil\log_2\gamma\rceil$ further probes and, as shown in the proof of \cref{thm:correct}, ends with a strict improvement of the incumbent; so the searches occur exactly on the $I$ improving pairs. This gives $e+I\lceil\log_2\gamma\rceil$, and $I\le e\le|P_1||P_2|$ gives the first bound in \cref{eq:N}. Since $P_3\ne\emptyset$, $|P_1|+|P_2|\le mn-1$, and for a fixed sum the product is largest when the two sizes differ by at most one, which gives $N$. After the relabelling, $|P_1|+|P_2|\le\frac23mn$, so $|P_1||P_2|\le\big(\frac{|P_1|+|P_2|}{2}\big)^2\le(mn)^2/9$. With integer times, the incumbent values are distinct integers, the first is at most $\sum_k\max_{P_k}t_{ij}$ and all are at least $\OPTNI\ge\sum_k\tau_k$ (\cref{lem:idle}), which gives the bound on $I$.

A probe evaluated by \cref{lem:probe} is an uncapacitated minimum cost flow problem on $V=m+n$ nodes and $E=mn$ arcs, which Orlin's algorithm solves in $O(V\log V\,(E+V\log V))$ time~\cite{Orlin1993}. Evaluated by \cref{lem:scc}, it is one integral maximum flow on $m+n+2$ nodes and at most $mn+m+n$ arcs, which takes $O(VE)=O((m+n)mn)$ time~\cite{Orlin2013}, followed by a linear-time component search~\cite{Tarjan1972}. Sorting takes $O(mn\log(mn))$ time. Multiplying the probe count, at most $(mn)^2(1+\lceil\log_2(mn)\rceil)/9$ since $\gamma\le mn$, by the cost of a probe gives the last statement.
\end{proof}

The quantity $N$ is the bound stated in~\cite{Jain2021} for the number of auxiliary time minimization problems solved by its algorithm, as the product $\alpha\beta$ of the numbers of distinct Phase-1 and Phase-2 times. A direct count gives a slightly larger bound. The recorded Phase-1 times strictly decrease, so at most $\alpha+1$ proficient triads are generated. Each costs one problem in Step~3 and at most $\beta$ problems in Step~5, because the Phase-2 values tested while generating one proficient triad are distinct and lie below the Phase-2 time of its first triad. With the two problems of Step~1, at most $(\alpha+1)(\beta+1)+2\le N+mn+2$ problems are solved, since $\alpha+\beta\le mn-1$. The product $\alpha\beta$ itself can be exceeded: for $a=(2,1)$, $b=(1,2)$, $t=\bigl(\begin{smallmatrix}3&2\\6&2\end{smallmatrix}\bigr)$, $P_1=\{(1,1)\}$, $P_2=\{(1,2),(2,1)\}$ and $P_3=\{(2,2)\}$, the algorithm solves six problems, four of them after Step~1, whatever optimal solutions are chosen, whereas $\alpha\beta=N=2$. \Cref{alg:tw} therefore has the same order of auxiliary computations, up to the factor $1+\lceil\log_2\gamma\rceil\le1+\lceil\log_2(mn)\rceil$, while enforcing non-idleness. The refined count $e+I\lceil\log_2\gamma\rceil$ explains the behavior observed in \cref{sec:experiments}: improvements are rare, so the number of probes is close to the number of evaluated pairs.

\begin{remark}[Effect of the cap]\label{rem:cap}
Consider the variant of \cref{alg:tw} without the cap, in which every pair that passes the margin test of line~10 is probed at $\theta^\gamma$ and, if $\Pi(w_1,w_2,\theta^\gamma)$ holds, $\theta^\star(w_1,w_2)$ is located by binary search on $\theta^1,\dots,\theta^\gamma$ and the incumbent is replaced only if $t_{w_1}+t_{w_2}+\theta^\star(w_1,w_2)<z^{\mathrm{inc}}$. Both versions generate, prune and evaluate the same pairs and produce the same sequence of incumbent values: by induction over the pairs, a pair updates the incumbent in either version exactly when $t_{w_1}+t_{w_2}+\theta^\star(w_1,w_2)<z^{\mathrm{inc}}$, and to the same value. The cap removes every probe that locates $\theta^\star$ for a pair that cannot improve the incumbent. On a pair that does improve it, the binary search runs over $h\le\gamma$ values instead of $\gamma$, which can take more or fewer probes. The bound $|P_1|\,|P_2|\,(1+\lceil\log_2\gamma\rceil)$ of \cref{thm:complexity} holds for both versions; the refined bound $e+I\lceil\log_2\gamma\rceil$ holds for the version with the cap. In the experiments of \cref{sec:capexp} the cap reduced the total number of probes in every instance.
\end{remark}

\subsection{Minimum phase times under non-idleness}

The algorithm of~\cite{Jain2021} starts from the minimum attainable time of each phase. In \PTP{} this minimum can be zero; under non-idleness the relevant quantity is
\begin{equation}\label{eq:TkminNI}
T^{\mathrm{NI}}_{k,\min}=\min\{T_k(X):X\in\SNI\},\qquad k=1,2,3,
\end{equation}
with $T^{\mathrm{NI}}_{k,\min}=+\infty$ if $\SNI=\emptyset$; in \cref{thm:tkmin} the minimum of an empty set is likewise $+\infty$.
By \cref{lem:idle}, $T^{\mathrm{NI}}_{k,\min}\ge\tau_k$. It can be computed with the same primitive as $\Pi$. For $k\in\{1,2,3\}$ let $\{a,b\}=\{1,2,3\}\setminus\{k\}$. For $w\in P_k$ and $w_a\in P_a$, say that $\Pi_k(w,w_a)$ holds if, after reserving one unit on $w$ and one on $w_a$, the transportation problem with the reduced margins and admissible links $\{(i,j)\in P_k:t_{ij}\le t_w\}\cup P_a\cup P_b$ has an integer feasible solution with positive Phase-$b$ flow.

\begin{theorem}\label{thm:tkmin}
Under \cref{ass:positive}, for each $k\in\{1,2,3\}$,
\[
T^{\mathrm{NI}}_{k,\min}=\min\{t_w:\ w\in P_k,\ \Pi_k(w,w_a)\text{ holds for some }w_a\in P_a\}.
\]
\end{theorem}

\begin{proof}
If $\Pi_k(w,w_a)$ holds with solution $X'$, then $X=X'+\ind_w+\ind_{w_a}\in\SNI$ and $T_k(X)=t_w$, exactly as in the proof of \cref{thm:twowitness}. Conversely, if $X\in\SNI$ attains $T^{\mathrm{NI}}_{k,\min}$, choose $w\in P_k$ with $x_w>0$ and $t_w=T_k(X)$, and any $w_a\in P_a$ with $x_{w_a}>0$; then $X-\ind_w-\ind_{w_a}$ shows that $\Pi_k(w,w_a)$ holds. The first part shows that $T^{\mathrm{NI}}_{k,\min}$ is at most the right-hand side, and the second that the right-hand side is at most $t_w=T^{\mathrm{NI}}_{k,\min}$. If $\SNI=\emptyset$, the first part shows that $\Pi_k(w,w_a)$ holds for no pair, so both sides are $+\infty$.
\end{proof}

\begin{corollary}\label{cor:empty}
$\SNI\neq\emptyset$ if and only if $T^{\mathrm{NI}}_{k,\min}$ is finite for some (equivalently, every) $k\in\{1,2,3\}$. Moreover, $\OPTNI\ge\sum_{k=1}^3T^{\mathrm{NI}}_{k,\min}$.
\end{corollary}

\begin{proof}
Each $T^{\mathrm{NI}}_{k,\min}$ is finite exactly when $\SNI\ne\emptyset$. If $X^*$ is optimal for~\cref{eq:NI}, then $T_k(X^*)\ge T^{\mathrm{NI}}_{k,\min}$ for each $k$, and summing over $k$ gives $\OPTNI=\sum_{k=1}^3T_k(X^*)\ge\sum_{k=1}^3T^{\mathrm{NI}}_{k,\min}$. If $\SNI=\emptyset$, both sides are $+\infty$.
\end{proof}

\begin{remark}\label{rem:tkprune}
\Cref{cor:empty} suggests replacing $\tau_k$ by $T^{\mathrm{NI}}_{k,\min}$ in the pruning bounds of \cref{alg:tw} and skipping witnesses with $t_w<T^{\mathrm{NI}}_{k,\min}$. This is correct. If $\Pi(w_1,w_2,\theta)$ holds, the solution $X$ of the first part of the proof of \cref{thm:twowitness} lies in $\SNI$ and has $T_1(X)=t_{w_1}$, $T_2(X)=t_{w_2}$ and $T_3(X)\le\theta$, so $t_{w_1}\ge T^{\mathrm{NI}}_{1,\min}$, $t_{w_2}\ge T^{\mathrm{NI}}_{2,\min}$ and $\theta\ge T^{\mathrm{NI}}_{3,\min}$. Witnesses below these bounds therefore never satisfy $\Pi$, every candidate value of a pair is at least $t_{w_1}+t_{w_2}+T^{\mathrm{NI}}_{3,\min}$, and the safety arguments of the proof of \cref{thm:correct} apply with $\tau_k$ replaced by $T^{\mathrm{NI}}_{k,\min}$; the cap index still satisfies $h\ge1$ because $\theta^1=\tau_3\le T^{\mathrm{NI}}_{3,\min}$. Computing $T^{\mathrm{NI}}_{k,\min}$ by \cref{thm:tkmin}, however, costs probes of its own. In the experiments of \cref{sec:capexp} the bound exceeded $\tau_k$ too rarely for this to pay: it increased the number of probes in 98\% of the instances (median increase 0.3\%).
\end{remark}

The bound $T^{\mathrm{NI}}_{k,\min}\ge\tau_k$ can be strict, so computing $T^{\mathrm{NI}}_{k,\min}$ is not a matter of reading off the smallest Phase-$k$ entry.

\begin{proposition}\label{prop:strict}
There are instances with $T^{\mathrm{NI}}_{k,\min}>\tau_k$.
\end{proposition}

\begin{proof}
Let $m=n=2$, $a=(1,2)$, $b=(1,2)$, $t=\begin{pmatrix}2&1\\4&3\end{pmatrix}$, $P_1=\{(1,2)\}$, $P_2=\{(1,1),(2,2)\}$ and $P_3=\{(2,1)\}$. Every $X\in\Shat$ has the form $x_{11}=s$, $x_{12}=1-s$, $x_{21}=1-s$, $x_{22}=1+s$ with $s\in\{0,1\}$. For $s=1$ the only Phase-1 link carries nothing, so $\SNI$ consists of the single solution with $s=0$, whose Phase-2 time is $t_{22}=3$. Hence $T^{\mathrm{NI}}_{2,\min}=3>2=\tau_2$. Using the fast Phase-2 link $(1,1)$ would exhaust $a_1=1$, and origin~1 is the only origin with a Phase-1 link.
\end{proof}

\subsection{Phase-flow bounds do not decide feasibility}

A natural shortcut would be to decide whether $\SNI\neq\emptyset$ from the largest flows that subsets of phases can carry. For every non-empty $K\subseteq\{1,2,3\}$, the condition
\begin{equation}\label{eq:flowbound}
\max_{X\in\Shat}\ \sum_{k\in K}F_k(X)\ \ge\ |K|
\end{equation}
is clearly necessary for $\SNI\neq\emptyset$. Each of these seven conditions can be checked with one minimum cost flow computation. The following example shows that they are not sufficient, which is why \cref{alg:tw} forces two phases explicitly through witnesses.

\begin{proposition}\label{prop:counter}
There is an instance satisfying~\cref{eq:flowbound} for every non-empty $K\subseteq\{1,2,3\}$ for which $\SNI=\emptyset$.
\end{proposition}

\begin{proof}
Let $m=n=2$, $a=(1,2)$, $b=(1,2)$, $P_1=\{(2,2)\}$, $P_2=\{(1,2),(2,1)\}$ and $P_3=\{(1,1)\}$ (the times play no role). As in the proof of \cref{prop:strict}, $\Shat=\{X_0,X_1\}$ with $x_{11}=s$, $x_{12}=x_{21}=1-s$, $x_{22}=1+s$ for $s=0,1$. Then $(F_1,F_2,F_3)=(1,2,0)$ for $X_0$ and $(2,0,1)$ for $X_1$. The left-hand side of~\cref{eq:flowbound} is $2,2,1$ for $K=\{1\},\{2\},\{3\}$, it is $3,3,2$ for $K=\{1,2\},\{1,3\},\{2,3\}$, and it is $3$ for $K=\{1,2,3\}$, so every condition holds. Yet neither solution is non-idle.
\end{proof}

The example also shows that the linear relaxation of $\SNI$ can be empty while each phase individually can be made busy: $F_3\ge1$ forces $s=1$, while $F_2\ge1$ forces $s\le\tfrac12$.

\subsection{Dependence on the number of phases}\label{sec:pphases}

The witness construction does not depend on there being three phases. For an integer $p\ge2$, let the link set be partitioned into non-empty sets $P_1,\dots,P_p$, define $F_k$ and $T_k$ by \cref{eq:FT} for $k=1,\dots,p$, and call
\begin{equation}\label{eq:pNI}
\OPTNI^{(p)}=\min\Big\{\textstyle\sum_{k=1}^{p}T_k(X)\ :\ X\in\Shat,\ F_k(X)\ge1,\ k=1,\dots,p\Big\}
\end{equation}
the $p$-phase problem with non-idle phases ($p$-PhTP-NI). For witnesses $w_k\in P_k$, $k<p$, and $\theta\in\Tset_p$, the reduced margins, the admissible link set $L(w_1,\dots,w_{p-1},\theta)$ and the condition $\Pi(w_1,\dots,w_{p-1},\theta)$ are defined as in \cref{def:aux}, with $p-1$ witnesses and with Phase-$p$ in the role of Phase-3.

\begin{proposition}\label{prop:fixedp}
Under \cref{ass:positive}, let $p\ge2$, relabel the phases so that $|P_p|\ge|P_k|$ for every $k$, and let $\gamma_p=|\Tset_p|$. Then
\[
\begin{aligned}
\OPTNI^{(p)}=\min\Big\{&\textstyle\sum_{k<p}t_{w_k}+\theta\ :\ w_k\in P_k\ (k<p),\ \theta\in\Tset_p,\\
&\Pi(w_1,\dots,w_{p-1},\theta)\text{ holds}\Big\}.
\end{aligned}
\]
Enumerating witness tuples and locating the smallest admissible $\theta$ by binary search requires at most
\begin{equation}\label{eq:pbound}
\prod_{k<p}|P_k|\,(1+\lceil\log_2\gamma_p\rceil)\ \le\ \Big(\frac{mn}{p}\Big)^{p-1}(1+\lceil\log_2\gamma_p\rceil)
\end{equation}
probes, so $p$-PhTP-NI is solvable in strongly polynomial time for every fixed $p$.
\end{proposition}

\begin{proof}
By \cref{rem:symmetry}, applied to $p$ phases, the relabelling does not change $\OPTNI^{(p)}$. Let witnesses $w_k\in P_k$, $k<p$, and $\theta\in\Tset_p$ satisfy $\Pi(w_1,\dots,w_{p-1},\theta)$ with solution $X'$, and let $X=X'+\sum_{k<p}\ind_{w_k}$. The witnesses lie in distinct phases, hence are distinct links, and adding one unit on each restores the original margins, so $X\in\Shat$. Every positive entry of $X$ lies in $L(w_1,\dots,w_{p-1},\theta)$ and $x_{w_k}\ge1$, so $T_k(X)=t_{w_k}$ for $k<p$, $T_p(X)\le\theta$, and $F_p(X)=F_p(X')\ge1$; thus $X$ is non-idle and $\OPTNI^{(p)}\le\sum_{k<p}t_{w_k}+\theta$. Conversely, if no non-idle solution exists, no tuple satisfies $\Pi$ by the first part, and both sides of the formula are $+\infty$. Otherwise, since $\Shat$ is finite, \cref{eq:pNI} has an optimal solution $X^*$. For $k<p$ choose $w_k\in P_k$ with $x^*_{w_k}>0$ and $t_{w_k}=T_k(X^*)$, and put $\theta=T_p(X^*)\in\Tset_p$. Then $X^*-\sum_{k<p}\ind_{w_k}$ is a non-negative integer matrix with the reduced margins, supported on $L(w_1,\dots,w_{p-1},\theta)$, with Phase-$p$ flow $F_p(X^*)\ge1$, so $\Pi(w_1,\dots,w_{p-1},\theta)$ holds with value $\sum_kT_k(X^*)=\OPTNI^{(p)}$. This proves the formula.

\Cref{lem:scc,lem:monotone} concern a single phase and hold with Phase-$p$ in place of Phase-3. For a fixed witness tuple, one probe at $\max\Tset_p$ decides whether some $\theta$ is admissible, and, by monotonicity, a binary search on the $\gamma_p$ values of $\Tset_p$ locates the smallest admissible one with at most $\lceil\log_2\gamma_p\rceil$ further probes, each one maximum flow and one strongly connected components computation. This gives at most $\prod_{k<p}|P_k|\,(1+\lceil\log_2\gamma_p\rceil)$ probes. Since $|P_p|\ge|P_k|$ for every $k$, $|P_p|\ge mn/p$, so $\sum_{k<p}|P_k|\le mn(p-1)/p$, and the inequality of arithmetic and geometric means gives $\prod_{k<p}|P_k|\le(mn/p)^{p-1}$. For fixed $p$ this number of probes is polynomial in $m$ and $n$, and each probe is strongly polynomial.
\end{proof}

For $p=3$, \cref{eq:pbound} reduces to the bound $(mn)^2(1+\lceil\log_2\gamma\rceil)/9$ of \cref{thm:complexity}. For $p=2$ a single witness suffices.

\begin{corollary}\label{cor:twophase}
Let $p=2$, relabel the phases so that $|P_1|\le|P_2|$, and let $\gamma_2=|\Tset_2|$. Then $\OPTNI^{(2)}$ is computed with at most $|P_1|\,(1+\lceil\log_2\gamma_2\rceil)$ probes, each consisting of one maximum flow and one strongly connected components computation on a network with $m+n+2$ nodes. No transportation costs are involved.
\end{corollary}
\begin{proof}
For $p=2$, the method of \cref{prop:fixedp} enumerates the witnesses $w_1\in P_1$ and, for each of them, locates the smallest admissible Phase-2 time with at most $1+\lceil\log_2\gamma_2\rceil$ probes, hence with at most $|P_1|\,(1+\lceil\log_2\gamma_2\rceil)$ probes in all. By \cref{lem:scc} with Phase-2 in the role of Phase-3, each probe is one maximum flow on the network with $m+n+2$ nodes, which yields a feasible $X^0$ or shows that none exists, and one strongly connected components computation on $D(X^0)$. Neither step uses transportation costs.
\end{proof}

\begin{remark}\label{rem:mincost}
If each link also has a unit cost $c_{ij}$, the minimum cost $C^*$ over the non-idle solutions with $T_1+T_2=\OPTNI^{(2)}$ equals
\[
\min\big\{c_{w_1}+c_{w_2}+\mathrm{MC}(w_1,w_2)\ :\ w_1\in P_1,\ w_2\in P_2,\ t_{w_1}+t_{w_2}=\OPTNI^{(2)}\big\},
\]
where $\mathrm{MC}(w_1,w_2)$ is the minimum cost of the transportation problem with the margins reduced at $w_1$ and $w_2$ on the links $\{(i,j)\in P_k:t_{ij}\le t_{w_k}\}$, $k=1,2$, and $+\infty$ if it has no feasible solution. Choosing as witnesses positive links of a non-idle solution of cost $C^*$ and sum of times $\OPTNI^{(2)}$ that attain $T_1$ and $T_2$ gives $t_{w_1}+t_{w_2}=\OPTNI^{(2)}$, and removing one unit from each witness leaves a feasible solution of the reduced problem, so the minimum is at most $C^*$. Conversely, adding the two witnesses to a feasible solution of such a reduced problem gives a non-idle solution with $T_k=t_{w_k}$, $k=1,2$, and hence with sum of times $\OPTNI^{(2)}$, so the minimum is at least $C^*$. At most $|P_1||P_2|$ minimum cost flow computations are needed. The pair $(\OPTNI^{(2)},C^*)$ is the non-dominated point of minimum sum of times of the bi-objective model of~\cite{Kaur2026arxiv} (\cref{sec:othermodels}).
\end{remark}

The exponent $p-1$ in \cref{eq:pbound} grows with the number of phases. We now show that, when $p$ is part of the input, no algorithm runs in polynomial time unless P${}={}$NP. Recall that an instance of three-dimensional matching (3DM) consists of disjoint sets $W$, $X$, $Y$ with $|W|=|X|=|Y|=q$ and a set $\mathcal{T}\subseteq W\times X\times Y$ of triples, and asks whether some $q$ triples of $\mathcal{T}$ cover every element of $W\cup X\cup Y$ exactly once. 3DM is NP-complete~\cite{Karp1972,GareyJohnson1979}.

\paragraph{Decision problems and encoding}
We use the terminology of~\cite{GareyJohnson1979}. A set of $q$ triples of $\mathcal{T}$ that covers every element of $W\cup X\cup Y$ exactly once is called a \emph{perfect matching}. An instance of 3DM is encoded by listing the elements of $W$, $X$ and $Y$ and the triples of $\mathcal{T}$, so its length is at least $q$ and at least $|\mathcal{T}|$. An instance $\mathcal{I}$ of $p$-PhTP-NI consists of $m$, $n$ and $p$, the availabilities $a_i$ and the requirements $b_j$, the times $t_{ij}$, which are positive rational numbers given by a numerator and a denominator in lowest terms, and, for every link $(i,j)$, the index of the phase that contains it; all integers are encoded in binary. Since the times of all $mn$ links are listed, the length of $\mathcal{I}$ satisfies $\mathrm{Length}(\mathcal{I})\ge mn$, and $mn\ge p$ because the $p$ phases are non-empty and disjoint. For the $p$-phase problem we write $\tau_k=\min\{t_{ij}:(i,j)\in P_k\}$ for $k=1,\dots,p$, and $\SNI$ for the set $\{X\in\Shat:F_k(X)\ge1,\ k=1,\dots,p\}$ of non-idle solutions, which is the feasible set of~\cref{eq:pNI}. We consider two decision problems. The \emph{feasibility problem} asks, for an instance $\mathcal{I}$, whether $\SNI\neq\emptyset$. An instance of the \emph{decision version} of $p$-PhTP-NI consists of an instance of $p$-PhTP-NI together with an integer $K$, encoded in binary; it is again denoted by $\mathcal{I}$, and the question is whether some $X\in\SNI$ satisfies $\sum_{k=1}^{p}T_k(X)\le K$. We write $\mathrm{Max}(\mathcal{I})$ for the largest integer that occurs among the $a_i$, the $b_j$, the numerators and denominators of the $t_{ij}$ and, in the decision version, $|K|$. A decision problem is \emph{strongly NP-complete} (NP-complete in the strong sense) if it belongs to NP and there is a polynomial $\pi$ such that its restriction to the instances with $\mathrm{Max}(\mathcal{I})\le\pi(\mathrm{Length}(\mathcal{I}))$ is NP-complete~\cite[Section~4.2]{GareyJohnson1979}. Our reduction from 3DM needs a restricted form.

\begin{lemma}\label{lem:3dm}
3DM remains NP-complete when no two triples share the same pair of $W$- and $X$-elements.
\end{lemma}

\begin{proof}
The restricted problem is a special case of 3DM and therefore belongs to NP. We show that it is NP-hard by a polynomial-time reduction from 3DM.

Let $(W,X,Y,\mathcal{T})$ be an instance of 3DM, let $\mathcal{T}_x$ be the set of triples containing $x\in X$, and let $d_x=|\mathcal{T}_x|$. If some $d_x=0$, the element $x$ cannot be covered, so the instance has no perfect matching, and we output the fixed instance with $q=1$ and no triples; it satisfies the restriction, vacuously, and has no perfect matching. So assume $d_x\ge1$ for every $x$. Construct $X'=\{x_\tau:\tau\in\mathcal{T}\}$, with one new element per triple; add new elements $\alpha_{x,i}$ to $W$ and $\beta_{x,i}$ to $Y$ for every $x\in X$ and $i=1,\dots,d_x-1$, giving $W'$ and $Y'$; and let $\mathcal{T}'$ consist of a \emph{main} triple $(w,x_\tau,y)$ for each $\tau=(w,x,y)\in\mathcal{T}$ and a \emph{filler} triple $(\alpha_{x,i},x_\tau,\beta_{x,i})$ for every $x\in X$, $\tau\in\mathcal{T}_x$ and $i\le d_x-1$. Every triple contains exactly one element of $X$, so $\sum_xd_x=|\mathcal{T}|$. Hence $|W'|=|X'|=|Y'|=q+\sum_x(d_x-1)=|\mathcal{T}|$, because $|X|=q$, and
\[
|\mathcal{T}'|=|\mathcal{T}|+\sum_{x\in X}d_x(d_x-1)=\sum_{x\in X}d_x^2\le\Big(\sum_{x\in X}d_x\Big)^2=|\mathcal{T}|^2 .
\]
The new instance thus has $3|\mathcal{T}|$ elements and at most $|\mathcal{T}|^2$ triples, and it is computed in time polynomial in $q+|\mathcal{T}|$, hence in the length of the given instance. Pairs are unique: in a main pair $(w,x_\tau)$ the element $x_\tau$ determines $\tau$, in a filler pair $(\alpha_{x,i},x_\tau)$ the element $\alpha_{x,i}$ determines $\beta_{x,i}$, and main and filler pairs have different kinds of first element.

It remains to show that $(W,X,Y,\mathcal{T})$ has a perfect matching if and only if $(W',X',Y',\mathcal{T}')$ has one. First let $M$ be a perfect matching of the original instance, and for $x\in X$ let $\tau_x$ be the triple of $M$ that contains $x$; then $M=\{\tau_x:x\in X\}$, and $\tau_x$ is the only triple of $\mathcal{T}_x$ in $M$. For each $x$ number the $d_x-1$ triples of $\mathcal{T}_x\setminus\{\tau_x\}$ as $\tau_{x,1},\dots,\tau_{x,d_x-1}$, and let $M'$ consist of the main triples of all $\tau_x$ and of the filler triples $(\alpha_{x,i},x_{\tau_{x,i}},\beta_{x,i})$, $x\in X$, $1\le i\le d_x-1$; these are triples of $\mathcal{T}'$ because $\tau_{x,i}\in\mathcal{T}_x$. An element of $W$ or of $Y$ lies only in main triples, and it lies in exactly one triple of $M$, hence in exactly one triple of $M'$. Each $\alpha_{x,i}$ and each $\beta_{x,i}$ lies in exactly one triple of $M'$, the filler triple with index $(x,i)$. Finally let $\tau\in\mathcal{T}$ and let $x$ be its element of $X$. If $\tau=\tau_x$, then $x_\tau$ lies in the main triple of $\tau$ and in no filler triple of $M'$; if $\tau=\tau_{x,i}$, then $x_\tau$ lies in the filler triple with index $(x,i)$, in no other filler triple, and in no main triple of $M'$, because $\tau\notin M$. So $M'$ covers $W'\cup X'\cup Y'$ exactly once.

Conversely, let $M'\subseteq\mathcal{T}'$ be a perfect matching of $(W',X',Y',\mathcal{T}')$, fix $x\in X$, and let $X'_x=\{x_\tau:\tau\in\mathcal{T}_x\}$, so that $|X'_x|=d_x$. An element $x_\tau\in X'_x$ occurs only in the main triple of $\tau$ and in the filler triples $(\alpha_{x,i},x_\tau,\beta_{x,i})$, $1\le i\le d_x-1$. Each $\alpha_{x,i}$ occurs only in filler triples of $x$, so $M'$ contains exactly one filler triple for each $i=1,\dots,d_x-1$; since $M'$ covers every element of $X'$ exactly once, these $d_x-1$ triples cover $d_x-1$ distinct elements of $X'_x$, and no other triple of $M'$ covers them. Exactly one element of $X'_x$ remains; call it $x_{\tau_x}$. It must be covered by $M'$, and not by a filler triple: every filler triple that contains it contains some $\alpha_{x,i}$, and the only triple of $M'$ containing $\alpha_{x,i}$ covers a different element of $X'_x$. Hence the main triple of $\tau_x$ lies in $M'$, and no other main triple of a $\tau\in\mathcal{T}_x$ lies in $M'$. Since every triple of $\mathcal{T}$ belongs to exactly one $\mathcal{T}_x$, the main triples in $M'$ are exactly those of $M=\{\tau_x:x\in X\}$. The triples $\tau_x$ are distinct, because they contain distinct elements of $X$, so $M$ is a set of $q$ triples in which each $x\in X$ lies in exactly one triple. Every element of $W$ and of $Y$ occurs in $\mathcal{T}'$ only in main triples and is covered exactly once by $M'$, hence by exactly one triple of $M$. Therefore $M$ is a perfect matching of $(W,X,Y,\mathcal{T})$.

The map is therefore a polynomial-time reduction from 3DM to its restriction, and the lemma follows from the NP-completeness of 3DM~\cite{Karp1972,GareyJohnson1979}.
\end{proof}

\begin{theorem}\label{thm:hard}
Let $p$ be part of the input.
\begin{itemize}
\item[(a)] Deciding whether $\SNI\neq\emptyset$ is strongly NP-complete, even when $m=n=p$, $a_i=b_j=1$ and $t_{ij}=1$ for all links.
\item[(b)] Given an integer $K$, deciding whether some $X\in\Shat$ with $F_k(X)\ge1$ for $k=1,\dots,p$ satisfies $\sum_{k=1}^pT_k(X)\le K$ is strongly NP-complete, even when $m=n=p$, $a_i=b_j=1$, $t_{ij}\in\{1,2\}$ and $K=p$.
\end{itemize}
\end{theorem}

\begin{proof}
Part~(a) concerns the feasibility problem and part~(b) the decision version of $p$-PhTP-NI.

\emph{Membership in NP.} Let $\mathcal{I}$ be an instance of either problem. Every $X\in\Shat$ satisfies $0\le x_{ij}\le\min\{a_i,b_j\}\le\mathrm{Max}(\mathcal{I})$, so $X$ is written with at most $mn\,(1+\log_2\mathrm{Max}(\mathcal{I}))$ bits. This is at most $\mathrm{Length}(\mathcal{I})\,(1+\mathrm{Length}(\mathcal{I}))$, because $\mathrm{Length}(\mathcal{I})\ge mn$ and because the largest integer of $\mathcal{I}$ is written in binary, so that $\log_2\mathrm{Max}(\mathcal{I})\le\mathrm{Length}(\mathcal{I})$. Given $X$, the row and column sums, the flows $F_k(X)$, the times $T_k(X)$ and their sum are obtained with $O(mn)$ additions and comparisons of integers and rational numbers whose encodings have length polynomial in $\mathrm{Length}(\mathcal{I})$. Hence a non-idle $X$ (with $\sum_kT_k(X)\le K$ in part~(b)) is a certificate verifiable in polynomial time, and both problems, and their restrictions considered below, belong to NP.

\emph{Three facts.}
(i) \Cref{lem:idle} and its proof remain valid for $k=1,\dots,p$: if $F_k(X)\ge1$, some link of $P_k$ is positive, so $T_k(X)\ge\tau_k>0$, and if $F_k(X)=0$ then $T_k(X)=0$.
(ii) If $a_i=b_j=1$ for all $i$ and $j$, then $m=n$ and $\Shat$ is the set of permutation matrices of order $n$: a row of non-negative integers with sum~1 has one entry equal to~1 and zeros elsewhere, and likewise every column. Each $X\in\Shat$ then has exactly $n$ positive links, one in each row and one in each column.
(iii) If, moreover, $p=n$, then $X\in\Shat$ is non-idle if and only if every phase contains exactly one positive link of $X$: the phases partition the link set, so $\sum_{k=1}^{p}F_k(X)=n=p$, and $F_k(X)\ge1$ for all $k$ forces $F_k(X)=1$ for all $k$.

\emph{Two fixed instances.} Both have $m=n=p=2$, $a=b=(1,1)$, all times equal to~1 and, in part~(b), $K=2$; they satisfy every restriction of~(a) and of~(b). In $\mathcal{I}^-$ the phases are $P_1=\{(1,1),(2,2)\}$ and $P_2=\{(1,2),(2,1)\}$. Each of the two permutation matrices of order~2 lies in a single phase, so $\SNI=\emptyset$ by~(ii), and $\mathcal{I}^-$ is a no-instance of both problems. In $\mathcal{I}^+$ the phases are $P_1=\{(1,1)\}$ and $P_2=\{(1,2),(2,1),(2,2)\}$. The identity matrix has $F_1=F_2=1$ and $T_1=T_2=1$, so $\mathcal{I}^+$ is a yes-instance of both problems.

Both reductions start from 3DM with unique pairs, which is NP-complete by \cref{lem:3dm}. Let $(W,X,Y,\mathcal{T})$ with $|W|=|X|=|Y|=q$ be an instance with unique pairs. If some $y\in Y$ lies in no triple, the instance has no perfect matching, and both reductions output $\mathcal{I}^-$. Assume from now on that every $y\in Y$ lies in some triple.

(a) If $q=0$, the empty set of triples is a perfect matching, and we output $\mathcal{I}^+$. Let $q\ge1$. Add a row $w^*$, a column $x^*$ and a phase $y^*$: the rows are $W\cup\{w^*\}$, the columns $X\cup\{x^*\}$, all margins are~1, and the $p=q+1$ phases are indexed by $Y\cup\{y^*\}$. Link $(w,x)$ is placed in phase $y$ if $(w,x,y)\in\mathcal{T}$, which is well defined because pairs are unique, and every other link, in particular every link of row $w^*$ and of column $x^*$, is placed in phase $y^*$. All times are~1. Every phase is non-empty: phase $y\in Y$ contains the link of a triple that contains $y$, and phase $y^*$ contains $(w^*,x^*)$. Moreover $m=n=p=q+1\ge2$, and the instance, which has $(q+1)^2$ links, is computed in time polynomial in $q+|\mathcal{T}|$.

If $M$ is a perfect matching, let $x_{wx}=1$ for every $(w,x,y)\in M$, $x_{w^*x^*}=1$, and let all other entries be~0. Each $w\in W$ and each $x\in X$ lies in exactly one triple of $M$, so this is a permutation matrix. Phase $y\in Y$ contains the link of the triple of $M$ that contains $y$, and phase $y^*$ contains $(w^*,x^*)$; so the matrix is non-idle. Conversely, let $X\in\SNI$. By (ii) and~(iii) it has $q+1$ positive links, one in each row, one in each column and one in each phase. All links of row $w^*$ and all links of column $x^*$ lie in phase $y^*$, which contains only one positive link, so the positive links of row $w^*$ and of column $x^*$ coincide, that is, $x_{w^*x^*}=1$. The other $q$ positive links therefore lie in $W\times X$, one in each row $w\in W$ and one in each column $x\in X$, and they occupy the $q$ phases $y\in Y$, one each. A link of $W\times X$ in a phase $y\in Y$ is the link of a triple $(w,x,y)\in\mathcal{T}$, so these links correspond to $q$ triples with distinct $w$ (rows), distinct $x$ (columns) and distinct $y$ (phases), that is, to a perfect matching.

The map is thus a polynomial-time reduction from 3DM with unique pairs to the feasibility problem restricted to the instances with $m=n=p$, $a_i=b_j=1$ and $t_{ij}=1$ for all links. This restricted problem is therefore NP-hard and, being in NP, NP-complete. All its instances have $\mathrm{Max}(\mathcal{I})=1$, so the same map is a reduction to the feasibility problem restricted to the instances with $\mathrm{Max}(\mathcal{I})\le1$, which is therefore NP-complete as well. Hence the feasibility problem is strongly NP-complete, and so is its restriction to the instances with $m=n=p$, $a_i=b_j=1$ and $t_{ij}=1$ for all links, all of whose instances have $\mathrm{Max}(\mathcal{I})=1$.

(b) If $q\le1$, the instance has a perfect matching: the empty set if $q=0$, and, if $q=1$, the triple that contains the element of $Y$, which is the only possible triple. We output $\mathcal{I}^+$. Let $q\ge2$. Take rows $W$, columns $X$, all margins equal to~1, and one phase for each $y\in Y$, so that $p=q$. For each triple $(w,x,y)$, place link $(w,x)$ in phase $y$ with time~1; this is well defined because pairs are unique. Give every other link time~2 and place it in one fixed phase. Every phase is non-empty, because every $y\in Y$ lies in some triple, $m=n=p=q\ge2$ and $K=q=p$. The instance has $q^2$ links and is computed in time polynomial in $q+|\mathcal{T}|$.

If the 3DM instance has a perfect matching $M$, let $x_{wx}=1$ for every $(w,x,y)\in M$ and let all other entries be~0. This is a permutation matrix that uses only links of time~1, exactly one in each phase, so every $T_k=1$, the solution is non-idle, and $\sum_kT_k=q$. Conversely, let $X$ be non-idle with $\sum_kT_k(X)\le q$. By~(i), $T_k(X)\ge\tau_k\ge1$ for every $k$, because every time is~1 or~2; as there are $q$ phases, every $T_k(X)=1$. A positive link of $X$ in phase $k$ has time at most $T_k(X)=1$, so every positive link of $X$ has time~1 and is therefore the link of a triple. By (ii) and~(iii), $X$ has $q$ positive links, one in each row, one in each column and one in each phase, and the link $(w,x)$ of a triple $(w,x,y)$ lies in phase $y$. The corresponding triples are therefore disjoint in $W$ (rows), $X$ (columns) and $Y$ (phases), so they form a perfect matching.

The map is thus a polynomial-time reduction from 3DM with unique pairs to the decision version of $p$-PhTP-NI restricted to the instances with $m=n=p$, $a_i=b_j=1$, $t_{ij}\in\{1,2\}$ and $K=p$, which is therefore NP-hard and, being in NP, NP-complete. Every such instance has $\mathrm{Max}(\mathcal{I})\le\max\{2,p\}=p\le\mathrm{Length}(\mathcal{I})$, so the same map is a reduction to the decision version restricted to the instances with $\mathrm{Max}(\mathcal{I})\le\mathrm{Length}(\mathcal{I})$, which is therefore NP-complete as well. Hence the decision version of $p$-PhTP-NI is strongly NP-complete, and so is its restriction to the instances with $m=n=p$, $a_i=b_j=1$, $t_{ij}\in\{1,2\}$ and $K=p$, all of whose instances have $\mathrm{Max}(\mathcal{I})\le\mathrm{Length}(\mathcal{I})$.
\end{proof}

\begin{remark}\label{rem:bfroma}
In the instances of part~(a) of \cref{thm:hard} all times equal~1. Every phase of a non-idle solution $X$ contains a positive link, whose time is~1, so $T_k(X)=1$ for all $k$ and $\sum_{k=1}^{p}T_k(X)=p$. Adding the bound $K=p$ to such an instance therefore gives an instance of part~(b) with the same answer. Together with the membership in NP shown in the proof, part~(a) thus implies part~(b), even with $t_{ij}=1$ for all links.
\end{remark}

\Cref{prop:fixedp,thm:hard} together locate the problem: $p$-PhTP-NI is solvable in strongly polynomial time for every fixed $p$, whereas for $p$ part of the input even the existence of a non-idle solution is strongly NP-complete. In the instances of \cref{thm:hard}, $p=m=n$, so an algorithm running in time $O((mn)^c)$, with $c$ and the constant hidden in the $O$ independent of $p$, would decide an NP-complete problem in polynomial time. Hence, unless P${}={}$NP, no such algorithm exists, and no polynomial-time algorithm can even return a non-idle solution whenever one exists. This excludes running times that are polynomial in the input size uniformly in $p$, but not running times of the form $f(p)\cdot\mathrm{poly}(m,n)$ for a computable function $f$. For the feasibility problem such running times are attained (\cref{prop:fptfeas} below). We are not aware of an algorithm with such a running time for the optimization problem $p$-PhTP-NI, nor of a result that excludes one. In the reduction of part~(b), the trivial bound $\OPTNI^{(p)}\ge p$, which follows from \cref{lem:idle}, is attained exactly when a perfect matching exists. Exhaustive computational checks of \cref{lem:3dm}, \cref{thm:hard} and \cref{prop:fixedp} on random instances are reported in Section~S7 of the Supplementary Material.

The feasibility problem has a simple combinatorial form, which places it among known fixed-parameter tractable problems.

\begin{lemma}\label{lem:feasw}
For $p$ phases, $\SNI\neq\emptyset$ if and only if there are links $w_k\in P_k$, $k=1,\dots,p$, such that every row $i$ contains at most $a_i$ of them and every column $j$ at most $b_j$.
\end{lemma}

\begin{proof}
If $X\in\SNI$, choose a positive link $w_k$ of $P_k$ for every $k$. These links are distinct, and since $X$ is integral, $X-\sum_k\ind_{w_k}\ge0$, so row $i$ contains at most $a_i$ of them and column $j$ at most $b_j$. Conversely, given such links, the margins reduced by one unit at every $w_k$ are non-negative integers with equal totals $Q-p$, so the northwest corner rule gives an integer solution $X'$ of the transportation problem with these margins on all links, and $X'+\sum_k\ind_{w_k}\in\SNI$.
\end{proof}

\begin{proposition}\label{prop:fptfeas}
Let $p$ be part of the input. Whether $\SNI\neq\emptyset$ can be decided, and an element of $\SNI$ found when one exists, with $2^{O(p\log p)}\cdot(mn)^{O(1)}$ operations on integers of at most $\mathrm{Length}(\mathcal{I})$ bits, hence in time $2^{O(p\log p)}\cdot\mathrm{Length}(\mathcal{I})^{O(1)}$. In particular, the feasibility problem is fixed-parameter tractable with parameter~$p$.
\end{proposition}

\begin{proof}
\emph{Reduction.} Let $G$ be the bipartite graph whose vertices are the row copies $(i,s)$ with $0\le s<\min\{a_i,p\}$ and the column copies $(j,u)$ with $0\le u<\min\{b_j,p\}$, and which has, for every link $(i,j)$ and all such $s$ and $u$, an edge between $(i,s)$ and $(j,u)$ of color $k$, where $P_k$ is the phase of $(i,j)$. A matching of $p$ edges of distinct colors, a \emph{rainbow matching} of size $p$~\cite{GuptaRSZ2017}, uses every color once and hence one link $w_k$ of every phase; the edges at a row use distinct copies of it, so row $i$ contains at most $\min\{a_i,p\}\le a_i$ of the $w_k$, and likewise for columns. Conversely, the links $w_k$ of \cref{lem:feasw} lie at most $\min\{a_i,p\}$ times in row $i$ and at most $\min\{b_j,p\}$ times in column $j$, so distinct copies can be assigned to them. By \cref{lem:feasw}, $\SNI\neq\emptyset$ if and only if $G$ has a rainbow matching of size $p$. Since $p\le mn$, $G$ has at most $mnp^2\le(mn)^3$ edges, and it is built with $O(m+n)$ comparisons of binary numbers and $O((mn)^3)$ further operations.

\emph{Kernel.} As in~\cite[Section~4.1]{GuptaRSZ2017}, identify each edge with the $3$-set formed by its two end vertices and its color; a rainbow matching of size $p$ is then a family of $p$ pairwise disjoint sets. A \emph{sunflower} with $r$ petals is a family of $r$ sets whose pairwise intersections all equal one set, the core. By the sunflower lemma of Erd\H{o}s and Rado~\cite{ErdosRado1960}, a family of more than $3!\,(r-1)^3$ distinct $3$-sets contains a sunflower with $r$ petals, and one is found in polynomial time: choose greedily a maximal family $\mathcal{G}$ of pairwise disjoint sets; if $|\mathcal{G}|\ge r$, then $r$ of its sets form a sunflower with empty core; otherwise every set meets one of the at most $3(r-1)$ elements of the sets of $\mathcal{G}$, so some element $x$ lies in more than $2(r-1)^2$ sets, and a sunflower with $r$ petals among the distinct $2$-sets $S\setminus\{x\}$ with $x\in S$, found in the same way (with singletons at the last level), gives one among the original sets after adding $x$ back. Let $r=3p-1$ and let $S_0$ be a member of such a sunflower. Then $p$ pairwise disjoint sets exist if and only if they exist after $S_0$ is deleted~\cite[Reduction Rule~12 and Lemma~13]{GuptaRSZ2017}: in a solution containing $S_0$, the other $p-1$ sets avoid the core, and since the petals are pairwise disjoint they meet at most $3(p-1)$ of the other $3p-2$ members of the sunflower, so one of these can replace $S_0$. Deleting sets in this way while more than $6(3p-2)^3$ sets remain takes at most $mnp^2$ rounds, each polynomial in $mn$, and leaves a subfamily with at most $6(3p-2)^3$ sets and the same answer.

\emph{Search.} The subfamily has at most $\binom{6(3p-2)^3}{p}\le(6(3p-2)^3)^p=2^{O(p\log p)}$ subsets with $p$ members, and each is tested for pairwise disjointness with $O(p^2)$ comparisons of vertex names and colors, which are integers or pairs of integers below $mn$. A subset that passes consists of edges of $G$ and gives links $w_k$ as in \cref{lem:feasw}, and the proof of the lemma turns them into an element of $\SNI$ with the northwest corner rule, which uses $O(m+n)$ additions and comparisons of integers at most $Q$. If no subset passes, the subfamily has no $p$ pairwise disjoint sets. By the kernel step, neither has the family of all edges of $G$, so $G$ has no rainbow matching of size $p$, and $\SNI=\emptyset$ by the reduction. The reduction and the kernel take $(mn)^{O(1)}$ operations and the search $2^{O(p\log p)}$ operations, which gives the stated bound.
\end{proof}

The reduction in the proof shows that faster algorithms for rainbow matching apply as well; for instance, the randomized algorithm with one-sided error of~\cite[Theorem~2]{GuptaRSZ2017} solves the feasibility problem in time $2^p\cdot\mathrm{Length}(\mathcal{I})^{O(1)}$. \Cref{prop:fptfeas} is consistent with \cref{thm:hard}, whose instances have $p=m=n$.

\section{Numerical illustrations}\label{sec:examples}

\subsection{An illustrative dispatch scenario}\label{sec:dispatch}

The following instance is synthetic; it was constructed for illustration and is not taken from operational data. A distributor serves stores $S_1,\dots,S_4$ from warehouses $W_1,W_2,W_3$ in three dispatch waves that share one fleet (\cref{tab:dispatch}); a unit is one truckload, and $t_{ij}$ is the round-trip time of route $(i,j)$ in hours. The example rests on three explicit assumptions. (i)~Each route may be used in one wave only, because of receiving windows and vehicle restrictions; this exclusive assignment gives the three phases. (ii)~There are enough vehicles to run all loads of a wave in parallel: all loads of a wave leave together and the next wave leaves when the last vehicle has returned, so a wave lasts as long as its longest round trip. (iii)~Every wave must carry at least one load because of an external commitment, such as a crew roster or a carrier contract, not because using every wave is desirable in itself.

\begin{table}[htbp]
\centering
\caption{Synthetic dispatch instance: round-trip times in hours and truckloads. Wave-1 routes are underlined, wave-2 routes are in bold, and wave-3 routes are plain.}\label{tab:dispatch}
\begin{tabular}{lccccc}
\toprule
 & $S_1$ & $S_2$ & $S_3$ & $S_4$ & $a_i$\\
\midrule
$W_1$ & \underline{2} & \underline{3} & \textbf{5} & 7 & 4\\
$W_2$ & \underline{3} & \textbf{2} & 6 & \textbf{4} & 3\\
$W_3$ & \textbf{8} & 7 & \underline{3} & 2 & 5\\
\midrule
$b_j$ & 3 & 4 & 2 & 3 & 12\\
\bottomrule
\end{tabular}
\end{table}

Without the commitment, the unique optimal plan ships $W_1\to S_2$ (4 loads), $W_2\to S_1$ (3) and $W_3\to S_3$ (2) in wave~1 and $W_3\to S_4$ (3) in wave~3, with triad $(3,0,2)$ and $\OPT=5$ hours; wave~2 is idle. With it, $\OPTNI=7$ hours with triad $(3,2,2)$: one to three loads use the 2-hour wave-2 route $W_2\to S_2$, and the wave-1 deliveries of $W_1$ and $W_2$ to $S_1$ and $S_2$ are rebalanced, for instance $x_{11}=1$, $x_{12}=3$, $x_{21}=2$, $x_{22}=1$, $x_{33}=2$, $x_{34}=3$. The commitment costs 2 hours, 40\% of the unrestricted duration. In general, $\OPTNI-\OPT$ is the price of requiring every phase to be used, which a planner can weigh against the cost of an idle phase. These values were verified by complete enumeration of the 317 feasible plans and by an implementation of \cref{alg:tw}.

\subsection{A small instance}\label{sec:small}

Consider the $2\times3$ instance with $a=(2,5)$, $b=(3,2,2)$ and the times and phases in \cref{tab:small}. Here $P_1=\{(1,3),(2,1)\}$, $P_2=\{(1,2),(2,2)\}$ and $P_3=\{(1,1),(2,3)\}$, so $\Tset_3=\{5,9\}$ and $\tau_2=2$, $\tau_3=5$.

\begin{table}[htbp]
\centering
\caption{Data of the small instance. Phase-1 times are underlined, Phase-2 times are in bold, and Phase-3 times are plain.}\label{tab:small}
\begin{tabular}{lcccc}
\toprule
 & $D_1$ & $D_2$ & $D_3$ & $a_i$\\
\midrule
$S_1$ & 9 & \textbf{4} & \underline{7} & 2\\
$S_2$ & \underline{6} & \textbf{2} & 5 & 5\\
\midrule
$b_j$ & 3 & 2 & 2 & 7\\
\bottomrule
\end{tabular}
\end{table}

The \PTP{} optimum is $x_{13}=2$, $x_{21}=3$, $x_{22}=2$, with triad $(7,2,0)$ and value~9; Phase-3 is idle, so this solution is not admissible for \PNI. \Cref{tab:smallrun} lists every step of \cref{alg:tw}; it solves four probes. The pair $((2,1),(2,2))$ shows a probe that fails because the auxiliary problem is infeasible, and the pair $((1,3),(2,2))$ shows a probe that is feasible but cannot route any flow through Phase-3. Because of the incumbent cap, each of the last two evaluated pairs is settled by a single probe at $\theta=5$; without the cap they would also be probed at $\theta=9$, for six probes in total.

\begin{table}[htbp]
\centering
\caption{Run of \cref{alg:tw} on the small instance. ``Bound'' is $t_{w_1}+t_{w_2}+\tau_3$; ``Cap'' is $\theta^h$, the largest Phase-3 time below $z^{\mathrm{inc}}-t_{w_1}-t_{w_2}$. A probe passes when $\Pi$ holds, and the number in parentheses is the largest Phase-3 flow found by the probe.}\label{tab:smallrun}
\small
\begin{tabular}{lllp{3.6cm}cc}
\toprule
$(w_1,w_2)$ & Bound & Cap & Probes & $\theta^\star$ & Value (incumbent)\\
\midrule
$((2,1),(2,2))$ & 13 & 9 & $\theta=9$: pass (4);\newline $\theta=5$: infeasible & 9 & 17 (17)\\
$((2,1),(1,2))$ & $15<17$ & 5 & $\theta=5$: pass (2) & 5 & 15 (15)\\
$((1,3),(2,2))$ & $14<15$ & 5 & $\theta=5$: fails (0) & -- & -- (15)\\
$((1,3),(1,2))$ & $16\ge15$ & \multicolumn{4}{l}{inner loop stops (pruned)}\\
\bottomrule
\end{tabular}
\end{table}

The optimal probe solution for $((2,1),(1,2),5)$ is $x'_{12}=1$, $x'_{21}=2$, $x'_{23}=2$. Adding the two reserved units gives $x_{12}=2$, $x_{21}=3$, $x_{23}=2$ with triad $(6,4,5)$, so $\OPTNI=15$. Non-idleness raises the optimal value from 9 to 15, and the adjustment is not simply the addition of Phase-3 flow: using the Phase-3 link $(2,3)$ releases origin~1 from the Phase-1 link $(1,3)$, so $T_1$ falls from 7 to 6 while $T_2$ rises from 2 to 4.

\subsection{Example 1 of Jain et al.}\label{sec:ex1}

\Cref{tab:ex1data} reproduces Example~1 of~\cite{Jain2021}, with $m=6$, $n=8$ and $Q=80$. Here $|P_1|=16$, $|P_2|=15$, $|P_3|=17$, and $\Tset_3$ has $\gamma=13$ values with $\tau_3=3$.

\begin{table}[htbp]
\centering
\caption{Example~1 of~\cite{Jain2021}. Phase-1 times are underlined, Phase-2 times are in bold, and Phase-3 times are plain.}\label{tab:ex1data}
\begin{tabular}{lccccccccc}
\toprule
 & $D_1$ & $D_2$ & $D_3$ & $D_4$ & $D_5$ & $D_6$ & $D_7$ & $D_8$ & $a_i$\\
\midrule
$S_1$ & \underline{06} & \textbf{03} & \textbf{12} & \underline{21} & \textbf{15} & 32 & 08 & 42 & 15 \\
$S_2$ & 40 & \underline{04} & 29 & 03 & \underline{23} & \underline{17} & \textbf{03} & \underline{06} & 10 \\
$S_3$ & \textbf{12} & 21 & \textbf{02} & \textbf{07} & 08 & \textbf{11} & \underline{13} & 24 & 12 \\
$S_4$ & \underline{08} & 05 & \underline{09} & \underline{04} & \textbf{09} & \underline{15} & 37 & \textbf{09} & 6 \\
$S_5$ & 29 & \textbf{27} & \textbf{17} & 19 & \underline{08} & 14 & \textbf{15} & \underline{23} & 17 \\
$S_6$ & 11 & \underline{35} & 42 & 05 & \underline{27} & \textbf{32} & \underline{35} & \textbf{27} & 20 \\
\midrule
$b_j$ & 12 & 5 & 13 & 6 & 9 & 20 & 8 & 7 & 80 \\
\bottomrule
\end{tabular}
\end{table}

The optimum of \PTP{} reported in~\cite{Jain2021} is $(06,32,00)$ with value~38; one optimal allocation, with phase flows $(F_1,F_2,F_3)=(19,61,0)$, is given in Table~S3 of the Supplementary Material. Phase-3 carries no flow. Since $|P_3|\ge\max\{|P_1|,|P_2|\}$, \cref{alg:tw} enumerates witnesses in Phase-1 and Phase-2. Of the $16\cdot15=240$ witness pairs it generates 187; 13 of these are pruned by the inner rule, one in each of the 13 inner loops that end early, and the remaining 174 are evaluated. It solves 180 probes in total; without the incumbent cap the same 174 pairs would require 721 (\cref{rem:cap}). The incumbent improves three times: to~48 with $w_1=(2,2)$, $w_2=(3,3)$, $\theta^\star=42$; to~45 with $w_1=(2,2)$, $w_2=(5,2)$, $\theta^\star=14$; and to~41 with $w_1=(1,1)$, $w_2=(6,6)$, $\theta^\star=3$. For the last pair the cap admits only $\theta<45-6-32=7$, so the probes test $\theta=5$ and $\theta=3$, and $\Pi$ holds at both. Each of the other 171 evaluated pairs is settled by a single probe that fails at the cap. The outer pruning rule never applies, because $t_{w_1}+\tau_2+\tau_3=t_{w_1}+5\le40<41$ for every Phase-1 link. One optimal allocation is shown in \cref{tab:ex1ni}. Optimal allocations are not unique: the implementation used in \cref{sec:experiments} returns a different one, with phase flows $(F_1,F_2,F_3)=(16,58,6)$ and the same triad.

\begin{table}[htbp]
\centering
\caption{An optimal allocation of \PNI{} for Example~1: triad $(06,32,03)$, value 41, phase flows $(F_1,F_2,F_3)=(7,67,6)$.}\label{tab:ex1ni}
\begin{tabular}{lccccccccc}
\toprule
 & $D_1$ & $D_2$ & $D_3$ & $D_4$ & $D_5$ & $D_6$ & $D_7$ & $D_8$ & $a_i$\\
\midrule
$S_1$ & \underline{7} & \textbf{5} & -- & -- & \textbf{3} & -- & -- & -- & 15 \\
$S_2$ & -- & -- & -- & 6 & -- & -- & \textbf{4} & -- & 10 \\
$S_3$ & \textbf{5} & -- & -- & -- & -- & \textbf{7} & -- & -- & 12 \\
$S_4$ & -- & -- & -- & -- & \textbf{6} & -- & -- & -- & 6 \\
$S_5$ & -- & -- & \textbf{13} & -- & -- & -- & \textbf{4} & -- & 17 \\
$S_6$ & -- & -- & -- & -- & -- & \textbf{13} & -- & \textbf{7} & 20 \\
\midrule
$b_j$ & 12 & 5 & 13 & 6 & 9 & 20 & 8 & 7 & 80 \\
\bottomrule
\end{tabular}
\end{table}

Thus $\OPTNI=41$ for Example~1. The Phase-1 and Phase-2 times are the same as in the \PTP{} optimum; non-idleness is met by sending 6 of the 80 units over the Phase-3 link $(2,4)$, whose time $3=\tau_3$ is the smallest in Phase-3. The cost of non-idleness is 3 time units, an increase of 7.9\%.

\paragraph{Filtering the proficient triads of~\cite{Jain2021}}
The supplementary material of~\cite{Jain2021} lists all proficient triads of Example~1 generated with respect to Phase-1, Phase-2 and Phase-3 (their Tables~10, 11 and~12). A tempting way to obtain the non-idle optimum is to discard the idle triads from such a list and keep the best remaining one. This gives 43 from the Phase-1 list (triad $(27,02,14)$) and 43 from the Phase-2 list, but 41 from the Phase-3 list (triad $(06,32,03)$, their row $l=9$). The result depends on the phase with respect to which the triads are generated, and two of the three choices miss the optimum. Moreover, the Phase-3 list contains the idle triad $(00,27,14)$, which also has value 41, so the value~41 is attained both by an admissible and by an inadmissible triad. Proficient triads are optimal for subproblems that allow idle phases, so a non-idle optimum need not appear among them. \Cref{sec:filterexp} shows that such failures are common in random instances drawn as described there.

\subsection{Example 2 of Jain et al.}\label{sec:ex2}

Example~2 of~\cite{Jain2021} has $m=n=10$ and $Q=143$. The published table marks Phase-1 links by underlining and Phase-2 links by encircling. We transcribed the partition from it; the transcribed data are given in Table~S4 of the Supplementary Material. As a check, the distinct times of each phase in our transcription coincide with the three ordered lists printed in~\cite{Jain2021}: $10<17<29<35<48<58$ for Phase-1, $12<21<31<39<45<61$ for Phase-2, and $9<19<28<43<60$ for Phase-3. The reported optimum of \PTP, $(10,21,09)$ with value~40, is reproduced by an exact enumeration of the threshold triads.

Here $|P_1|=36$ and $|P_2|=|P_3|=32$, so, by \cref{rem:symmetry}, \cref{alg:tw} enumerates witnesses in Phase-2 and Phase-3 and decides Phase-1 by the probe. It generates 100 of the $32\cdot32=1024$ witness pairs, prunes 8, evaluates 92, and solves 97 probes; the outer pruning rule ends the enumeration. The incumbent improves from 79 to 48 to 40, and the last improvement uses the Phase-2 witness $(1,1)$ with time 21 and the Phase-3 witness $(2,4)$ with time 9, for which the probe gives Phase-1 time 10. Hence $\OPTNI=40$ with triad $(10,21,09)$; the optimal solution found has phase flows $(F_1,F_2,F_3)=(55,39,49)$. In agreement with \cref{prop:relation}, the non-idleness constraint is inactive. The two published examples thus show both behaviors: the constraint binds in Example~1 and is inactive in Example~2.

\section{Computational experiments}\label{sec:experiments}

\subsection{Environment and implementation}

\Cref{alg:tw} was coded in Python, as in~\cite{Kaur2025matcom}, and run on one core of an Intel Xeon processor at 2.10~GHz with 7~GB of memory under Linux (Python 3.13, SciPy 1.18.1~\cite{Virtanen2020}). Probes were evaluated by \cref{lem:scc}, with SciPy's implementation of Dinic's maximum flow algorithm and its strongly connected components routine. CPU times are process times in seconds and exclude instance generation. Because the hardware differs, our times are not directly comparable with those of~\cite{Kaur2025matcom}.

Except in the solver comparison, whose times come from separate runs (\cref{sec:mipexp}; Section~S1 of the Supplementary Material), all CPU times reported below come from a single process pinned to one core, which re-ran every instance after the experiments; each re-run reproduced the logged probe count and optimal value. During this re-run the other core was busy with solver runs; in a separate test such concurrent load slowed a solve by about 6\%. Identical cells of different experiments differed in total CPU time by at most 2.4\%, and in the repeated runs of the solver comparison, on instances with a median time of at least 1~s, single runs deviated from their median by at most 21.8\% (at most 3.8\% at the median); shorter runs vary more. Section~S1 of the Supplementary Material gives the details. Probe counts are the exact measure of work.

\subsection{Test instances}

The design follows~\cite{Kaur2025matcom}. Times are independent uniform integers in $[1,r]$. In the first experiment the sizes are $10\times10$, $10\times15$, $20\times20$, $25\times20$, $30\times30$, $30\times35$, $40\times40$ and $r\in\{25,40,70,85,100\}$; in the second, for large instances, the sizes are $40\times45$ and $50\times50$ and $r\in\{50,300,500,1000,2000\}$. In both the links are split at random into three phases of equal size (rounded), the analogue of the 50\% split of~\cite{Kaur2025matcom}. In the third experiment $r=40$, the sizes are $10\times10$, $15\times15$, $20\times20$, $30\times30$, Phase-1 receives 5, 10, 25, 33, 50, 75, 90\% of the links, and the remaining links are split equally between Phase-2 and Phase-3. Since~\cite{Kaur2025matcom} does not state how availabilities and requirements are generated, we draw $a_i$ uniformly from $[10,50]$ and split $Q=\sum_ia_i$ at $n-1$ distinct uniform cut points of $\{1,\dots,Q-1\}$, so that every $b_j\ge1$. Each combination has 50 instances with recorded seeds, 3,650 instances in all, and each was also solved for \PTP{} to measure the effect of non-idleness. All instances have $\SNI\ne\emptyset$.

\subsection{Measures}

A witness pair is \emph{generated} when the loops of \cref{alg:tw} reach it, \emph{pruned} when a bound test eliminates it (including the pair that ends an inner loop), and \emph{evaluated} when at least one probe is solved for it; pairs with a negative reduced margin are counted separately (27,499 such pairs occurred). We report CPU time, probes, and the number $I$ of incumbent updates; means are given with standard deviations and 95\% $t$-intervals, proportions with 95\% Wilson intervals.

\subsection{Performance}

\Cref{tab:cpu} reports mean CPU times; per-size standard deviations and 95\% intervals are given in Table~S1 of the Supplementary Material. The largest instances, $50\times50$, took 26.4~s on average and at most 50.2~s. \Cref{tab:probes} gives the mean number of probes. Two features explain these counts. The pruning rules confine the enumeration to a shrinking part of the pair space: the fraction of pairs generated falls from 42\% at $10\times10$ to 7\% at $40\times40$. Improvements are rare (at most 18 in any instance), so by \cref{thm:complexity} almost every evaluated pair costs one probe: evaluated pairs account for 99.8\% of all probes (99.3\% per instance on average). Fitted by least squares over all instances of the first two experiments, the number of probes grows as $(mn)^\kappa$ with $\kappa=1.38$ (95\% interval $[1.37,1.40]$), well below the worst-case exponent 2 of \cref{thm:complexity}; CPU time grows with exponent 1.65 ($[1.64,1.66]$).

The width of the time range matters less than the size. Within each size, the Spearman correlation between $r$ and CPU time lies between 0.07 and 0.29. A wider range creates more distinct Phase-3 times and more improvements, so the binary searches are longer: for the first experiment the mean number of probes beyond one per evaluated pair is 10.2 for $r=25$ and 21.9 for $r=100$. Since this is small against the number of evaluated pairs, the running time is governed by the size of the pair space.

\begin{table}[htbp]
\centering
\caption{Mean CPU time in seconds of \cref{alg:tw} (50 instances per cell, equal phase sizes).}\label{tab:cpu}
\begin{tabular}{lrrrrr}
\toprule
Size & $[1,25]$ & $[1,40]$ & $[1,70]$ & $[1,85]$ & $[1,100]$\\
\midrule
$10\times10$ & 0.12 & 0.13 & 0.12 & 0.13 & 0.14\\
$10\times15$ & 0.22 & 0.23 & 0.23 & 0.27 & 0.25\\
$20\times20$ & 1.17 & 1.19 & 1.13 & 1.28 & 1.33\\
$25\times20$ & 1.63 & 1.67 & 1.88 & 1.92 & 1.88\\
$30\times30$ & 4.08 & 4.38 & 4.36 & 4.36 & 4.42\\
$30\times35$ & 5.32 & 5.39 & 5.59 & 5.52 & 6.03\\
$40\times40$ & 10.35 & 11.24 & 11.53 & 12.18 & 12.34\\
\midrule
Size & $[1,50]$ & $[1,300]$ & $[1,500]$ & $[1,1000]$ & $[1,2000]$\\
\midrule
$40\times45$ & 13.01 & 14.42 & 15.90 & 15.37 & 14.37\\
$50\times50$ & 23.70 & 25.21 & 27.62 & 27.33 & 27.89\\
\bottomrule
\end{tabular}
\end{table}

\begin{table}[htbp]
\centering
\caption{Witness pairs $|P_1||P_2|$ and mean number of probes (50 instances per cell).}\label{tab:probes}
\begin{tabular}{lrrrrrr}
\toprule
Size & Pairs & $[1,25]$ & $[1,40]$ & $[1,70]$ & $[1,85]$ & $[1,100]$\\
\midrule
$10\times10$ & 1,089 & 418 & 446 & 432 & 435 & 470\\
$10\times15$ & 2,500 & 747 & 762 & 796 & 850 & 811\\
$20\times20$ & 17,689 & 3,217 & 3,328 & 3,190 & 3,496 & 3,563\\
$25\times20$ & 27,722 & 4,375 & 4,476 & 4,954 & 4,986 & 4,979\\
$30\times30$ & 90,000 & 9,193 & 9,938 & 9,708 & 9,723 & 10,047\\
$30\times35$ & 122,500 & 11,273 & 11,397 & 11,734 & 11,487 & 12,562\\
$40\times40$ & 284,089 & 18,437 & 19,904 & 20,492 & 21,565 & 21,801\\
\midrule
Size & Pairs & $[1,50]$ & $[1,300]$ & $[1,500]$ & $[1,1000]$ & $[1,2000]$\\
\midrule
$40\times45$ & 360,000 & 22,511 & 24,400 & 26,240 & 25,189 & 23,697\\
$50\times50$ & 693,889 & 35,725 & 37,276 & 39,264 & 38,938 & 40,122\\
\bottomrule
\end{tabular}
\end{table}

\paragraph{Effect of the phase sizes}
\Cref{tab:share} varies the share of Phase-1. As in~\cite{Kaur2025matcom}, the running time is largest when the link set is split evenly: at every size the slowest share was 33\%. After the relabelling of line~1, \cref{alg:tw} enumerates the two smallest phases, whose product $|P_A||P_B|$ is largest, about $(mn)^2/9$, when the three phases have equal size; when Phase-1 is small it is one of the enumerated phases, and when it is large the other two phases are small.

\begin{table}[htbp]
\centering
\caption{Mean CPU time in seconds for different shares of Phase-1 links ($r=40$; Phase-2 and Phase-3 share the rest equally).}\label{tab:share}
\begin{tabular}{lrrrrrrr}
\toprule
Size & 5\% & 10\% & 25\% & 33\% & 50\% & 75\% & 90\%\\
\midrule
$10\times10$ & 0.029 & 0.049 & 0.093 & 0.126 & 0.057 & 0.014 & 0.007\\
$15\times15$ & 0.080 & 0.132 & 0.355 & 0.479 & 0.217 & 0.040 & 0.012\\
$20\times20$ & 0.191 & 0.360 & 0.937 & 1.171 & 0.491 & 0.069 & 0.021\\
$30\times30$ & 0.548 & 1.087 & 3.061 & 4.339 & 1.820 & 0.256 & 0.053\\
\bottomrule
\end{tabular}
\end{table}

\subsection{Effect of the incumbent cap}\label{sec:capexp}

To isolate the cap we solved 160 further instances (sizes 8, 12, 16 and 20; times in $[1,40]$ or $[10,20]$; availabilities in $[10,50]$ or $[1,5]$; 10 instances per combination) with and without it, counting every probe. By \cref{rem:cap} the two versions generate, prune and evaluate the same pairs, which the runs confirmed, and they returned the same optimal values. The cap reduced the number of probes in every instance, by 76.7\% on average (range 17.2--84.8\%); the median ratio of probes with and without it is 0.215 (5th--95th percentile 0.158--0.388; Wilcoxon signed-rank test, $p=5.2\times10^{-28}$; per-cell medians in Table~S2 of the Supplementary Material). Replacing $\tau_k$ by $T^{\mathrm{NI}}_{k,\min}$ in the pruning bounds did not pay (\cref{rem:tkprune}).

\subsection{How often non-idleness binds}

Over the first two experiments the constraint binds, that is $\OPTNI>\OPT$, in 926 of 2,250 instances (41.2\%, 95\% interval 39.1--43.2\%), with a mean relative increase of 4.1\% and a mean absolute increase of 1.60 time units (\cref{tab:bind}). In 14\% of the binding instances the computed optimum of \PTP{} uses a single phase. The relative cost falls as the range widens, from 8.0\% for $r=25$ to 0.5\% for $r=2000$. Within the first experiment, which crosses seven sizes with five time ranges, the binding rate does not depend significantly on the time range ($\chi^2=4.0$, 4 degrees of freedom, $p=0.40$) but depends on the size ($\chi^2=30.5$, 6 degrees of freedom, $p=3.1\times10^{-5}$): it is 55\% at $10\times10$, 45\% at $20\times20$ and 36--38\% from $30\times30$ to $40\times40$. In the second experiment it is 30\% at $40\times45$ and 38\% at $50\times50$. Beyond $30\times30$ the rate thus levels off between 30 and 40\% rather than continuing to decrease. The cost is small in relative terms because, with times starting at 1, every phase has a link that is nearly free; it grows when the smallest time of each phase is large, since $\OPTNI\ge\tau_1+\tau_2+\tau_3$ by \cref{lem:idle}.

\begin{table}[htbp]
\centering
\caption{Percentage of instances in which non-idleness binds / mean relative increase $(\OPTNI-\OPT)/\OPT$ in percent over those instances (50 instances per cell).}\label{tab:bind}
\begin{tabular}{lccccc}
\toprule
Size & $[1,25]$ & $[1,40]$ & $[1,70]$ & $[1,85]$ & $[1,100]$\\
\midrule
$10\times10$ & 62 / 7.3 & 44 / 5.6 & 54 / 4.7 & 60 / 3.0 & 56 / 3.3\\
$10\times15$ & 54 / 6.4 & 46 / 5.1 & 40 / 4.1 & 54 / 3.7 & 54 / 3.4\\
$20\times20$ & 46 / 7.6 & 52 / 4.3 & 46 / 3.8 & 40 / 2.8 & 42 / 2.3\\
$25\times20$ & 40 / 7.1 & 34 / 5.2 & 32 / 3.1 & 52 / 2.5 & 42 / 2.5\\
$30\times30$ & 36 / 8.2 & 36 / 5.0 & 36 / 2.9 & 40 / 2.7 & 34 / 2.6\\
$30\times35$ & 44 / 8.9 & 42 / 6.3 & 40 / 3.5 & 34 / 3.2 & 30 / 2.4\\
$40\times40$ & 42 / 11.3 & 42 / 6.8 & 30 / 4.0 & 38 / 2.9 & 38 / 3.1\\
\midrule
Size & $[1,50]$ & $[1,300]$ & $[1,500]$ & $[1,1000]$ & $[1,2000]$\\
\midrule
$40\times45$ & 34 / 5.3 & 22 / 1.1 & 34 / 0.8 & 26 / 0.7 & 34 / 0.6\\
$50\times50$ & 30 / 6.5 & 44 / 1.1 & 36 / 0.8 & 46 / 0.5 & 34 / 0.5\\
\bottomrule
\end{tabular}
\end{table}

\subsection{Correctness checks}\label{sec:correct}

On 5,000 further random instances of sizes $2\times2$, $2\times3$, $3\times2$, $3\times3$, $2\times4$, $4\times2$ and $3\times4$, availabilities in $\{1,2,3\}$, times in $\{1,\dots,9\}$ and random non-empty phases, 537 of which have $\SNI=\emptyset$, \cref{alg:tw} agreed with complete enumeration of $\Shat$ in every case, including the detection of infeasibility; the one-sided 95\% upper confidence bound on the probability of a disagreement is 0.06\%. This comparison is the primary check. For larger instances, where enumeration is impossible, we solved the two mixed integer formulations of \cref{sec:mipexp} with HiGHS 1.15.1~\cite{Huangfu2018} on 90 instances drawn from the first two experiments, 10 for each of the nine sizes. All 540 solver runs (two formulations, three runs each) proved optimality and returned the value found by \cref{alg:tw}.\footnote{Before the timing runs, both formulations were also checked against \cref{alg:tw} on 3,300 further random instances with $m,n\le8$, five of which have no non-idle solution. In that check, with the settings of \cref{sec:mipexp}, HiGHS 1.15.1 twice reported as optimal a value that was not optimal, once for each formulation: 24 where the optimum is 19, and 14 where it is 11. Both optima were established independently of \cref{alg:tw} and of the solver by an exact enumeration over phase thresholds, and for each saved model a feasible solution with the optimal value was checked against all its constraints. With presolve disabled, both runs returned the optimum. When the two saved models were solved again on another computer, the second discrepancy (14 instead of 11) recurred with HiGHS 1.13.0 to 1.15.1, but not with 1.11.0 or 1.12.0. The first (24 instead of 19) did not recur with the default options. With the option \texttt{random\_seed=15}, however, HiGHS 1.13.0 to 1.15.1 reported a different non-optimal value, 28, as optimal for the same saved model, whose optimum is 19; with presolve disabled they returned 19. We have reported this to the HiGHS developers as issue~\#3433 (\url{https://github.com/ERGO-Code/HiGHS/issues/3433})~\cite{HiGHSissue3433}. Which random seeds trigger the failure depends on the platform and on the HiGHS version; the runs, model files and logs are given in Section~S5 of the Supplementary Material. Agreement with the solver is therefore supporting evidence; the comparison with complete enumeration is the primary check.}

\subsection{Comparison with a mixed integer programming solver}\label{sec:mipexp}

Because \PNI{} has compact mixed integer formulations, a general-purpose solver is the natural alternative to \cref{alg:tw}. We compared \cref{alg:tw} with two formulations. The direct formulation is
\begin{equation}\label{eq:mip}
\begin{aligned}
\min\ & T_1+T_2+T_3\\
\text{s.t. } & X\in\Shat,\quad x_{ij}\le\min\{a_i,b_j\}\,y_{ij},\quad T_k\ge t_{ij}\,y_{ij}\ \ ((i,j)\in P_k),\\
& \textstyle\sum_{(i,j)\in P_k}x_{ij}\ge1,\quad y_{ij}\in\{0,1\},\qquad k=1,2,3,
\end{aligned}
\end{equation}
a big-$M$ model with $2mn+3$ variables, in which $y_{ij}=1$ marks a used link. Each constraint $T_k\ge t_{ij}y_{ij}$ involves a single link, so the linear relaxation, which can set $y_{ij}=x_{ij}/\min\{a_i,b_j\}$, bounds each $T_k$ weakly. A tighter formulation uses one binary variable per phase and distinct time instead of one per link. For phase $k$ let $d_{k,1}<\dots<d_{k,\gamma_k}$ be the distinct times of $P_k$ and $\ell(i,j)$ the level of link $(i,j)$, so that $t_{ij}=d_{k,\ell(i,j)}$; the binary variable $z_{k,l}$, $l=2,\dots,\gamma_k$, indicates $T_k\ge d_{k,l}$. The threshold formulation is
\begin{equation}\label{eq:mipz}
\begin{aligned}
\min\ & \textstyle\sum_{k=1}^{3}\big(d_{k,1}+\sum_{l=2}^{\gamma_k}(d_{k,l}-d_{k,l-1})\,z_{k,l}\big)\\
\text{s.t. } & X\in\Shat,\quad \textstyle\sum_{(i,j)\in P_k}x_{ij}\ge1\quad(k=1,2,3),\\
& x_{ij}\le\min\{a_i,b_j\}\,z_{k,\ell(i,j)}\quad((i,j)\in P_k,\ \ell(i,j)\ge2),\\
& z_{k,l}\le z_{k,l-1}\quad(3\le l\le\gamma_k),\qquad z_{k,l}\in\{0,1\}.
\end{aligned}
\end{equation}
For a non-idle $X$, setting $z_{k,l}=1$ exactly when $T_k(X)\ge d_{k,l}$ gives a feasible point whose objective equals $\sum_kT_k(X)$, because $T_k(X)\ge d_{k,1}$ by \cref{lem:idle}. Conversely, in a feasible $(X,z)$ a positive link of level $l\ge2$ forces $z_{k,l}=1$ and, by the monotonicity constraints, $z_{k,l'}=1$ for $2\le l'\le l$, so the objective is at least $\sum_kT_k(X)$. Hence the optimal value of~\eqref{eq:mipz} is $\OPTNI$. Model~\eqref{eq:mipz} has $mn$ integer flow variables and $\sum_k(\gamma_k-1)\le mn-3$ binary variables, and links of equal time in a phase share one binary variable.

We timed \cref{alg:tw} and HiGHS 1.15.1~\cite{Huangfu2018} with both models on the 90 instances of \cref{sec:correct}. HiGHS was called through its Python interface highspy 1.15.1 and restricted to one thread, with relative gap zero, a time limit of 600~s and all other options at their defaults; its time includes model construction. Each method was run three times per instance and the median CPU process time recorded, with the order of the three methods rotated across rounds and instances. On these 90 instances the times of \cref{alg:tw} in \cref{tab:mip} exceed those of the re-run behind \cref{tab:cpu} by 12.7\% at the median (interquartile range 9.5--17.7\%; Section~S1 of the Supplementary Material). The three methods in \cref{tab:mip} were measured under identical conditions, so the time ratios are not affected, but absolute times should be compared within a table. The node counts of HiGHS were identical across the three repetitions.

\Cref{alg:tw} was faster than model~\eqref{eq:mip} on 80 of the 90 instances, with a median time ratio of 2.8 (interquartile range 1.81--4.22; Wilcoxon signed-rank test, $p=3.8\times10^{-12}$; per-size medians 2.0--5.1). Model~\eqref{eq:mipz} is much stronger: the median number of branch-and-bound nodes was 1 to 3.5 for every size up to $40\times40$ and about 100 for $40\times45$ and $50\times50$, against medians of 54 to 6,661 nodes for model~\eqref{eq:mip}, and model~\eqref{eq:mipz} was faster than model~\eqref{eq:mip} on 70 instances (median time ratio 0.55, $p=4.0\times10^{-10}$). Against model~\eqref{eq:mipz}, \cref{alg:tw} was faster on 65 of the 90 instances, with a median time ratio of 1.5 (interquartile range 0.94--2.31; $p=9.1\times10^{-4}$). The per-size medians range from 0.8 at $40\times40$, where the solver was faster, to 3.0 at the two smallest sizes; for the 20 instances of sizes $40\times45$ and $50\times50$ the median is 2.0. Three factors bear on the comparison. Two favor the solver: HiGHS is compiled code, whereas our implementation is interpreted Python, and the solver needs no problem-specific code. The third is the formulation, whose effect is large: replacing~\eqref{eq:mip} by~\eqref{eq:mipz} roughly halved the solver's time at the median, and further strengthening or a commercial solver may narrow the gap more. Independently of measured times, \cref{thm:complexity} bounds the work of \cref{alg:tw} a priori, whereas the work of a branch-and-bound search admits no polynomial bound in general.

\begin{table}[htbp]
\centering
\caption{\Cref{alg:tw} against HiGHS (one thread) with the direct model~\eqref{eq:mip} and the threshold model~\eqref{eq:mipz}, 10 instances per size. Times in seconds; each method was run three times per instance and the median taken. Alg.: mean of the medians; \eqref{eq:mip}, \eqref{eq:mipz}: mean and maximum of the medians, including model construction; ratios: median over the instances of solver time divided by \cref{alg:tw} time. All solver runs reached proven optimality within the 600~s limit.}\label{tab:mip}
\small
\begin{tabular}{lrrrrrrrr}
\toprule
 & & & \multicolumn{2}{c}{Model \eqref{eq:mip}} & \multicolumn{2}{c}{Model \eqref{eq:mipz}} & \multicolumn{2}{c}{Ratio to Alg.}\\
\cmidrule(lr){4-5}\cmidrule(lr){6-7}\cmidrule(lr){8-9}
Size & $N$ & Alg. & Mean & Max & Mean & Max & \eqref{eq:mip} & \eqref{eq:mipz}\\
\midrule
$10\times10$ & 10 & 0.15 & 0.46 & 0.9 & 0.46 & 0.7 & 2.1 & 3.0\\
$10\times15$ & 10 & 0.32 & 0.72 & 1.8 & 0.91 & 1.5 & 2.0 & 3.0\\
$20\times20$ & 10 & 1.17 & 5.45 & 10.6 & 1.94 & 3.6 & 5.1 & 1.7\\
$25\times20$ & 10 & 1.79 & 6.32 & 10.9 & 2.52 & 4.6 & 3.5 & 1.2\\
$30\times30$ & 10 & 4.77 & 14.85 & 25.8 & 6.25 & 9.8 & 3.0 & 1.4\\
$30\times35$ & 10 & 6.44 & 15.66 & 27.0 & 7.71 & 17.1 & 2.8 & 1.2\\
$40\times40$ & 10 & 12.64 & 32.40 & 55.0 & 10.57 & 22.0 & 2.8 & 0.8\\
$40\times45$ & 10 & 17.92 & 70.07 & 184.6 & 33.32 & 51.4 & 3.7 & 2.1\\
$50\times50$ & 10 & 31.03 & 68.71 & 127.1 & 48.35 & 97.9 & 2.5 & 1.7\\
\bottomrule
\end{tabular}
\end{table}

\subsection{Filtering the proficient triads of Jain et al.}\label{sec:filterexp}

For 2,000 random instances with $m,n\in\{3,4\}$, availabilities in $\{1,\dots,4\}$, times in $\{1,\dots,12\}$, random phases in which Phase-1 receives about half of the links on average (the sampler is specified in Section~S4 of the Supplementary Material), and $\SNI\ne\emptyset$, we enumerated $\Shat$ and computed, for each phase $k$, the set $S^*_k$ of all triads that are proficient with respect to Phase-$k$, ties included. Every list that the algorithm of~\cite{Jain2021} can record when generating proficient triads with respect to Phase-$k$ is a subset of $S^*_k$, whatever its tie-breaking, by Theorem~4.5 of~\cite{Jain2021}, whose proof is the argument in the proof of \cref{prop:framework} with $\Shat$ in place of $\SNI$; so if discarding the idle triads of $S^*_k$ does not give $\OPTNI$, filtering fails for every implementation, and the frequency of this event is a lower bound on the failure rate. We also ran the algorithm of~\cite{Jain2021} itself, with exact auxiliary solutions and random tie-breaking. Filtering fails, whatever the implementation, with respect to at least one phase in 1,597 instances (79.8\%) and with respect to all three in 525 (26.2\%); for the algorithm as run, the numbers are 1,671 (83.5\%) and 965 (48.2\%) (\cref{tab:filter}). Among the 1,707 instances in which non-idleness binds, the lower bound for Phase-1 is 64.7\%. The rates of \cref{tab:filter} are the frequencies observed in this sample, and its intervals express their sampling error; the rates depend on the sample and on the sampling scheme. Twelve further samples of 2,000 instances from the same distribution gave Phase-1 lower bounds between 56.4\% and 61.4\%, all above the 55.2\% of this sample; pooled, they estimate the Phase-1 lower bound under this distribution at 59.2\% (95\% Wilson interval 58.5--59.8\%), which the Phase-1 interval of \cref{tab:filter} does not contain, so the Phase-1 rate of the table is lower than the estimated rate under the distribution. When each link instead receives an independent uniform phase label, eight samples of 2,000 instances gave Phase-1 lower bounds between 40.4\% and 42.9\% (41.8\% pooled). The seeds, the per-sample results and the rates for at least one phase are given in Section~S4.

\begin{table}[htbp]
\centering
\caption{Percentage of instances (95\% Wilson interval) in which discarding the idle triads from the proficient triads of~\cite{Jain2021} fails to give $\OPTNI$ (2,000 instances).}\label{tab:filter}
\small
\begin{tabular}{lccc}
\toprule
 & \multicolumn{3}{c}{Proficient triads with respect to}\\
\cmidrule(lr){2-4}
 & Phase-1 & Phase-2 & Phase-3\\
\midrule
All proficient triads $S^*_k$ & 55.2 (53.0--57.4) & 55.7 (53.5--57.9) & 55.5 (53.4--57.7)\\
\quad (lower bound) & & & \\
Recorded by~\cite{Jain2021} & 64.2 (62.0--66.2) & 67.8 (65.8--69.9) & 67.0 (65.0--69.1)\\
\quad (random ties) & & & \\
\bottomrule
\end{tabular}
\end{table}

\section{Discussion}\label{sec:discussion}

\subsection{Relation to the iterative framework of Jain et al.}

The algorithm of~\cite{Jain2021} generates proficient triads by solving restricted time minimization problems in which links are blocked by assigning them a prohibitively large time. Its auxiliary problems can make a link unusable, but they cannot require a link, or a phase, to carry flow. Non-idleness is a requirement of the second kind, and this is the main obstacle to adapting the iterative framework directly. The witnesses in \cref{alg:tw} meet the requirement by reserving one unit in each of two phases, which changes the margins rather than the link times, and the remaining phase is handled by maximizing its flow in the probe.

The convergence theory of~\cite{Jain2021} nevertheless carries over. Call a solution \emph{NI-restricted feasible} for an auxiliary problem of~\cite{Jain2021} if it is restricted feasible in the sense of their Definition~3.4 and belongs to $\SNI$, and call a triad \emph{NI-proficient} with respect to Phase-1 if $T_2+T_3$ is minimal among the non-idle solutions with Phase-1 time $T_1$. Consider Algorithm~3.1 of~\cite{Jain2021} with the following changes: (N1) every feasibility test and every auxiliary optimization of Steps~3--6 is taken over NI-restricted feasible solutions; (N2) the test of Step~2 compares $T1_{\mathrm{prof}}(l)$ with $T^{\mathrm{NI}}_{1,\min}$ instead of $T1_{\min}$, where $T^{\mathrm{NI}}_{1,\min}$ is computed in Step~1 by solving $(\mathrm{TMP}_1)^0$ over NI-restricted feasible solutions, or by \cref{thm:tkmin}; (N3) the algorithm stops at Step~1 if $\SNI=\emptyset$ (\cref{cor:empty}). Under the convention $T^{\mathrm{NI}}_{1,\min}=+\infty$ for $\SNI=\emptyset$, (N3) is implied by (N2). The Phase-2 floor used in Steps~4 and~5 may be $T2_{\min}$ of~\cite{Jain2021}, or $T^{\mathrm{NI}}_{2,\min}$, or any value of $\{0\}\cup\{t_{ij}:(i,j)\in P_2\}$ not exceeding $T^{\mathrm{NI}}_{2,\min}$. Both (N1) and (N2) are needed: without (N1) the algorithm may return a value below $\OPTNI$, or reach a step with no applicable branch; without (N2) it may reach Step~3 with no NI-restricted feasible solution, a case for which Algorithm~3.1 has no step.

\begin{proposition}\label{prop:framework}
With (N1)--(N3), and with ``feasible solution of 3-PhTP'' read as ``element of $\SNI$'' and ``proficient'' as ``NI-proficient'', Theorem~4.1, Corollary~4.2, Theorem~4.3 (with $T^{\mathrm{NI}}_{k,\min}$), Lemma~4.4, Theorem~4.5, Lemma~4.6, Theorem~4.7 and Corollary~4.8 of~\cite{Jain2021} remain valid. In particular the modified algorithm returns $\OPTNI$.
\end{proposition}

\begin{proof}
For $u\in(0,+\infty]$ let $F_u=\{X\in\SNI:T_1(X)<u\}$, and let $B(X)=T_2(X)+T_3(X)$.

\emph{The auxiliary problems.} In $(\mathrm{TMP}_3)^0(T_u)$ the time matrix is $t'_{ij}=t_{ij}$ on $P_3$ and $t'_{ij}=0$ on the unblocked links of $P_1\cup P_2$. Hence the objective $\max\{t'_{ij}:x_{ij}>0\}$ of every solution is $T_3(X)$. A non-idle $X$ is restricted feasible exactly when $T_1(X)<T_u$, because $T_1(X)$ is the time of a positive Phase-1 link and every positive Phase-1 link has time at most $T_1(X)$. So, under (N1), $(\mathrm{TMP}_3)^0(T_u)$ minimizes $T_3$ over $F_{T_u}$.

Let the current triad be that of $X^k$, put $B_k=B(X^k)$ and $c=t_2^{\beta_k-\delta}$. Then $(\mathrm{TMP}_3)^{k(\delta)}(T_u)$ minimizes $T_3$ over the $X\in F_{T_u}$ with $T_2(X)\le c$ and $T_3(X)<B_k-c$.

Fix $T_u$ with $F_{T_u}\neq\emptyset$, and consider any run of Steps~3--7 with any optimal solution chosen at each read-off. Let $X^0,\dots,X^s$ be its solutions, and let $c_k$ be the value of $c$ at which $X^{k+1}$ is found.

\emph{Lemma~4.4.}
\begin{itemize}
\item Relation~(3) holds because every $X^k$ lies in $F_{T_u}$.
\item From $T_2(X^{k+1})\le c_k<T_2(X^k)$ we get~(4).
\item From $T_3(X^{k+1})<B_k-c_k$ we get $B_{k+1}<B_k$, which is~(6).
\item For~(5): $T_3(X^1)\ge T_3(X^0)$, because $X^0$ minimizes $T_3$ over $F_{T_u}$. For $k\ge1$, if $T_3(X^{k+1})<T_3(X^k)$, then $T_2(X^{k+1})<T_2(X^k)\le c_{k-1}$ and $T_3(X^{k+1})<T_3(X^k)<B_{k-1}-c_{k-1}$. So $X^{k+1}$ is restricted feasible for the problem solved by $X^k$ and has a smaller objective, a contradiction.
\end{itemize}

\emph{Theorem~4.1.} An $X\in F_{T_u}$ with $T_2(X)=T_2(X^k)$ and $T_3(X)<T_3(X^k)$ would, by the same inequalities, be restricted feasible for the problem solved by $X^k$, with a smaller objective. This is impossible.

\emph{Theorem~4.5.} We show by induction on $k$: every $X\in F_{T_u}$ with $T_2(X)\ge T_2(X^k)$ satisfies $B(X)\ge B_k$.
\begin{itemize}
\item For $k=0$ this holds because $X^0$ minimizes $T_3$.
\item Let it hold for $k$, and let $X\in F_{T_u}$ have $T_2(X)\ge T_2(X^{k+1})$. If $T_2(X)\ge T_2(X^k)$, then $B(X)\ge B_k>B_{k+1}$. Otherwise suppose $B(X)<B_{k+1}$.
  \begin{itemize}
  \item If $T_2(X)>c_k$, then $T_2(X)$ is a Phase-2 time strictly between $c_k$ and $T_2(X^k)$, so the scan tested it before $c_k$. But $X$ is restricted feasible there, since $T_3(X)<B_k-T_2(X)$. This contradicts the choice of $c_k$ as the first restricted feasible value.
  \item If $T_2(X)\le c_k$, then $T_3(X)<T_3(X^{k+1})<B_k-c_k$. So $X$ is restricted feasible at $c_k$ with a smaller objective than $X^{k+1}$, a contradiction.
  \end{itemize}
\end{itemize}
Finally, let $X\in F_{T_u}$ with $T_2(X)<T_2(X^s)$, and suppose $B(X)<B_s$.
\begin{itemize}
\item If Step~7 was reached from Step~4, then $T_2(X^s)$ is the floor, which does not exceed $T_2(X)$. This is impossible.
\item Otherwise every Phase-2 time from the floor up to below $T_2(X^s)$ was tested and found infeasible. Among them is $T_2(X)$, at which $X$ is restricted feasible, a contradiction.
\end{itemize}
Hence $B_s=\min\{B(X):X\in F_{T_u}\}$, attained by $X^s$, and the last triad is NI-proficient with respect to Phase-1 (Remark~4.1 of~\cite{Jain2021}).

\emph{Well-definedness and termination.}
\begin{itemize}
\item Every Phase-2 time of a non-idle solution lies in the scanned set $\{t_2^0,t_2^1,\dots\}$ with $t_2^0=0$, and the floor is itself a value of this set that does not exceed the Phase-2 time of any non-idle solution. So Steps~4 and~5 always have an applicable branch.
\item The scanned values strictly decrease in each step, so Steps~3--7 terminate.
\item By (N2), Steps~2--7 continue exactly while $F_{T1_{\mathrm{prof}}(l)}\neq\emptyset$, and $F_\infty\ne\emptyset$ by (N3). Every recorded triad lies in $F_{T1_{\mathrm{prof}}(l)}$, so $T1_{\mathrm{prof}}$ strictly decreases, which is~(7), and the outer loop terminates.
\end{itemize}

\emph{Lemma~4.6\,(8) and Theorem~4.7.}
\begin{itemize}
\item The recorded triads are attained in $\SNI$.
\item For $l\ge1$, the $(l+1)$-st recorded solution lies in $F_{T1_{\mathrm{prof}}(l-1)}$, over which the $l$-th minimizes $B$. This gives~(8).
\item Let $X\in\SNI$, and let $l\ge1$ be the first index with $T1_{\mathrm{prof}}(l)\le T_1(X)$. It exists because the last recorded value is $T^{\mathrm{NI}}_{1,\min}$. Then $X\in F_{T1_{\mathrm{prof}}(l-1)}$, and the $l$-th recorded triad has $T1_{\mathrm{prof}}(l)+T2_{\mathrm{prof}}(l)+T3_{\mathrm{prof}}(l)\le T_1(X)+B(X)$.
\end{itemize}
So the minimum over the recorded triads is $\OPTNI$, which is Theorem~4.7. Step~8 returns this minimum, so the modified algorithm returns $\OPTNI$.

\emph{Theorem~4.3, Corollary~4.2 and Corollary~4.8.} Solved over NI-restricted feasible solutions, $(\mathrm{TMP}_k)^0$ returns $T^{\mathrm{NI}}_{k,\min}$; under (N1), $(\mathrm{TMP}_3)^0(T_u)$ returns $\min\{T_3(X):X\in F_{T_u}\}$. A non-idle minimizer of $T_1+T_2+T_3$ minimizes the sum of the two other phase times among the non-idle solutions with the same time in any one phase, which is Corollary~4.8.
\end{proof}

We also checked \cref{prop:framework} computationally. With every auxiliary problem solved exactly by enumeration, and the optimal solution used in each read-off of Steps~3 and~6 chosen at random among all optimal solutions, the modified algorithm returned $\OPTNI$ on all of 3,000 random instances with $m,n\le4$, and every relation of Section~4 of~\cite{Jain2021} held at every iteration. Without (N2) the algorithm reached a step with no applicable branch on 1,257 of the 3,000 instances (41.9\%).

The obstacle to using the modified algorithm is computational. Its auxiliary problems are time minimization problems with the three side constraints $F_k(X)\ge1$, so Remark~3.3 of~\cite{Jain2021}, which derives polynomiality from the auxiliary problems being standard TMTPs, no longer applies. The only exact method we know for these side-constrained problems is the witness construction of \cref{sec:algorithm}: with it, deciding whether such a problem has a non-idle solution with Phase-3 time at most $\theta$ costs, in the worst case, as many probes as there are witness pairs. One auxiliary problem of the modified algorithm can therefore cost, in the worst case, as much as \cref{alg:tw} applied to the whole instance, and the modified algorithm solves up to $(\alpha+1)(\beta+1)+2\le N+mn+2$ auxiliary problems. This is why we solve \PNI{} directly.

\subsection{Relation to other models}\label{sec:othermodels}

The restricted-flow model of~\cite{Kaur2025matcom} implies non-idleness whenever its prescribed phase quantities are positive, but fixes the split of the requirement across phases (\cref{rem:choices}). In \PNI{} the split is free, and the witnesses reserve only one unit per phase. The bi-objective model of~\cite{Kaur2026arxiv} has two phases and minimizes the sum of the phase times and the total cost, over solutions in which the flow of each phase is non-zero. Its algorithm starts from a cost-minimal solution and lowers the sum of times step by step, and by Theorem~2 of~\cite{Kaur2026arxiv} its last non-dominated point has the minimum sum of times over non-idle solutions. \Cref{cor:twophase} computes this sum directly, without costs, and \cref{rem:mincost} gives the cost of that point. The non-zero requirement is enforced in~\cite{Kaur2026arxiv} by a rule (its Remark~12): when the cost-minimal solution leaves a phase empty, the problem is solved again with some links blocked, and when the optimal solution of a sub-problem leaves a phase empty, the sub-problem is replaced by its immediate next minor sub-problem, the one with the next lower time bounds.

This rule does not guarantee the minimum. A sub-problem whose cheapest admissible solution leaves a phase empty may still contain non-idle solutions, possibly only at the corner of its time bounds; replacing it by a minor sub-problem discards them, and the stopping rule of Theorem~2 of~\cite{Kaur2026arxiv} can then certify a sum of times that is not minimal. The instance of \cref{tab:remark12} shows this. Row~2 ships its unit to one terminal and row~1 serves the other two, so $\Shat$ has three solutions: $X^a$ with $x_{21}=1$ has $(T_1,T_2)=(2,2)$ and cost~4; $X^b$ with $x_{22}=1$ has $(T_1,T_2)=(1,0)$ and cost~5 and leaves Phase-II empty; and $X^c$ with $x_{23}=1$ has $(T_1,T_2)=(2,1)$ and cost~8. The non-dominated points are $(4,4)$ and $(3,8)$, and the minimum sum of times is~3. The algorithm of~\cite{Kaur2026arxiv} starts from the cost-minimal solution $X^a$, with sum of times~4. Of its sub-problems restricted at the sum~4, the major ones bound the Phase-I and Phase-II times by $(2,1)$ and by $(1,2)$. In both, the cheapest admissible solution is $X^b$, which leaves Phase-II empty, so Remark~12 replaces each by an immediate next minor sub-problem. Every such sub-problem either bounds the Phase-I time by~1, which excludes the link $(1,2)$ used by $X^c$, or admits no Phase-II link. No non-idle solution with a smaller sum of times is found, and Theorem~2 of~\cite{Kaur2026arxiv} declares~4 the minimum. This holds for every choice of immediate next minor sub-problem and does not involve the treatment of the cost-minimal solution. Nor does it depend on how ties between optimal solutions are broken. Every sub-problem of~\cite{Kaur2026arxiv}, including the minor ones, minimizes the cost over the solutions that use no blocked link, where a Phase-I or Phase-II link is blocked when its time exceeds a bound for its phase (or reaches it, for a strict bound). Every link used by $X^b$ is a Phase-I link of time~1, and $X^c$ also uses such a link, so every sub-problem that admits $X^c$ also admits $X^b$, which is cheaper. Hence $X^c$ is never an optimal solution of a sub-problem, and since the cost-minimal solution $X^a$ is non-idle, the cost-minimal problem is not solved again.

\begin{table}[htbp]
\centering
\caption{A two-phase instance on which Remark~12 of~\cite{Kaur2026arxiv} certifies a non-minimal sum of times. Entries are time\,/\,cost; Phase-I links are underlined and Phase-II links are in bold.}\label{tab:remark12}
\begin{tabular}{lcccc}
\toprule
 & $D_1$ & $D_2$ & $D_3$ & $a_i$\\
\midrule
$S_1$ & \underline{1\,/\,3} & \underline{2\,/\,2} & \underline{1\,/\,1} & 2\\
$S_2$ & \textbf{2\,/\,1} & \underline{1\,/\,1} & \textbf{1\,/\,3} & 1\\
\midrule
$b_j$ & 1 & 1 & 1 & 3\\
\bottomrule
\end{tabular}
\end{table}

To estimate how often this occurs, we implemented the algorithm of~\cite{Kaur2026arxiv} step by step; on the numerical example of~\cite{Kaur2026arxiv} it reproduces the published non-dominated points exactly. We compared its output with complete enumeration on 10,000 random two-phase instances (sizes $2\times2$ to $3\times4$, total availability at most 8, times in $\{1,\dots,6\}$, random costs in $\{1,\dots,9\}$ for half of the instances and unit costs for the other half), 9,975 of which have non-idle solutions. Under the reading of Remark~12 that lowers the time bound of the phase carrying all the flow, the set of non-dominated points was wrong in 1,404 of the 9,975 instances (14.1\%, 95\% interval 13.4--14.8\%): in 825 the algorithm reported that no non-idle solution exists, in 310 it reported a minimum sum of times that is too large, and in 269 it missed a non-dominated point but ended at the correct minimum. Under the most permissive reading, which tries both immediate next minor sub-problems, 13.6\% of the sets were wrong. All failures occurred among the 2,520 instances in which Remark~12 was applied; when every sub-problem was instead solved exactly over non-idle solutions, the output was correct on all 10,000 instances. Because Remark~12 leaves open which links are blocked and which minor sub-problem is chosen, these rates depend on its reading. The instance of \cref{tab:remark12} does not, and Section~S6 of the Supplementary Material gives two $2\times2$ instances on which the rule fails for every choice of blocked links. These observations concern version~1 of the preprint~\cite{Kaur2026arxiv}. The preprint does not consider three phases.

The bilevel model of~\cite{ZhangXie2025} requires both groups of a partition of the origins to be non-empty, which is a condition on the partition rather than on the flow.

\subsection{Limitations and extensions}

The analysis assumes integer availabilities and requirements and positive times (\cref{ass:positive}). For $p$ phases, \cref{prop:fixedp} gives an algorithm whose exponent grows with $p$, and \cref{thm:hard} shows that, unless P${}={}$NP, no algorithm runs in time polynomial in the input size uniformly in $p$, since even deciding non-idle feasibility is strongly NP-complete for unbounded $p$; we are not aware of an algorithm for the optimization problem with running time $f(p)\cdot\mathrm{poly}(m,n)$, although the feasibility problem has one (\cref{prop:fptfeas}). Capacitated links, as in~\cite{XieLi2022,ZhangXie2024}, can be incorporated in the probe directly, since minimum cost flow handles arc capacities; the reserved units must then respect the capacities of the witness links. For operational use, the model also presumes that the units assigned to a phase can be shipped in parallel, so that fleet size is not a constraint, and that a phase lasts as long as its longest trip rather than a fixed clock interval (\cref{sec:dispatch}). A three-phase extension of the bi-objective model of~\cite{Kaur2026arxiv}, which would also minimize transportation cost, is a further direction.

\section{Conclusions}\label{sec:conclusions}

We studied the three-phase time minimization transportation problem with the additional requirement that every phase carries a positive quantity, while the split of the requirement across phases remains free. The requirement is not vacuous: the optimum of the main example of~\cite{Jain2021} leaves Phase-3 idle, and enforcing non-idleness raises its value from 38 to 41. For an operational reading, $\OPTNI-\OPT$ is the price of requiring every phase to be used (\cref{sec:dispatch}). We showed that the problem reduces to enumerating witness links in two phases and deciding the third by a single minimum cost flow computation, and that the optimal Phase-3 time for each pair can be found by binary search. The resulting algorithm is strongly polynomial and needs at most $N(1+\lceil\log_2\gamma\rceil)$ probes, where $N$ is the bound stated for the original algorithm in~\cite{Jain2021}; capping the binary search by the incumbent value ensures that a pair which cannot improve the incumbent costs a single probe. We also characterized the minimum attainable time of each phase, showed that it can exceed the smallest time entry of the phase, and showed that feasibility cannot be decided from phase-flow bounds alone.

With $p$ phases the problem is solvable in strongly polynomial time for each fixed $p$, and, when $p$ is part of the input, even deciding whether a non-idle solution exists is strongly NP-complete. Feasibility alone is fixed-parameter tractable in $p$ (\cref{prop:fptfeas}). We are not aware of an algorithm for the optimization problem with running time $f(p)\cdot\mathrm{poly}(m,n)$, nor of a result that excludes one. For two phases the method computes exactly the minimum-time point of the bi-objective model of~\cite{Kaur2026arxiv}, which the re-solution rule of that model can miss.

In computational experiments on random instances with up to $50\times50$ links, the number of probes grew as $(mn)^{1.38}$ on the 2,250 instances with equal phase sizes, far below the worst-case bound, and the incumbent cap removed 77\% of the probes on average on 160 further instances. Non-idleness raised the optimal value in 41\% of the instances with equal phase sizes, by 4.1\% on average. Discarding idle triads from the output of the original algorithm failed to give the non-idle optimum in more than half of the instances of every sample of small random instances drawn as in \cref{sec:filterexp}, and in more than 40\% of those of every sample in which each link received an independent phase label. An independent mixed integer programming solver confirmed the optimal values on 90 instances with up to $50\times50$ links; at the median it needed 2.8 times as long as the algorithm with a direct formulation and 1.5 times as long with a threshold formulation, which was the faster method on 25 of the 90 instances.

\section*{CRediT authorship contribution statement}
\textbf{Sailesh Kamalahasan:} Conceptualization, Methodology, Software, Investigation,
Data curation, Writing -- review \& editing, Visualization.
\textbf{Raghav Paladugu:} Conceptualization, Methodology, Software, Validation,
Formal analysis, Investigation, Data curation, Writing -- original draft, Visualization.
\textbf{Archana Khurana:} Resources, Supervision, Project administration.

\section*{Declaration of competing interest}
The authors declare that they have no known competing financial interests or personal
relationships that could have appeared to influence the work reported in this paper.

\section*{Funding}
This research did not receive any specific grant from funding agencies in the public,
commercial, or not-for-profit sectors.

\section*{Data and code availability}
The Supplementary Material consists of a document (Sections~S1--S8) and an archive. The archive contains scripts, recorded seeds and logs that reproduce the filtering experiment of \cref{sec:filterexp} and its two sensitivity analyses, the correctness check of \cref{sec:correct} on its archived instances with an independent implementation of \cref{alg:tw}, the count of auxiliary problems in the example of \cref{sec:theory} and the dispatch scenario of \cref{sec:dispatch}, together with the logged records of the HiGHS pre-check described in the footnote of \cref{sec:correct} and reruns of its two saved models. It also contains the original source code and recorded outputs of the computational study of \cref{sec:experiments}, including the solver comparison of \cref{sec:mipexp}, of the study of~\cite{Kaur2026arxiv} in \cref{sec:othermodels} and of the checks reported in Section~S7; these are historical records and were not rerun for the archive. For the computational check of \cref{prop:framework}, the archive contains the original source code and a log of an October 2026 reproduction using two seeded batches of 1,500 instances. The log records aggregate results, not individual instance outcomes; records of the original runs were not retained.

\bibliographystyle{elsarticle-num}
\bibliography{references}

\end{document}


\maketitle

This document supplements the article. It details the timing measurements (Section~S1), gives additional experimental tables (Section~S2) and the data of the two examples of Jain et al.\ (Section~S3), specifies the filtering study completely (Section~S4), documents the HiGHS observations and their reproduction (Section~S5), gives two further instances for Section~7.2 of the main text (Section~S6) and computational checks of results of Section~4.5 of the main text (Section~S7), and describes the accompanying archive of scripts, data and logs (Section~S8). Sections, tables, results and equations of the main text are cited with their numbers in the main text. The references of this document are numbered independently of those of the main text.

\section{Computing environment and timing}\label{S:timing}
\paragraph{Timing measurements (Section~6.1 of the main text)}

All CPU times reported in Section~6 of the main text, except those of the solver comparison of Section~6.8 of the main text, come from a single process pinned to one core, which re-ran every instance after the experiments; each re-run reproduced the logged probe count and optimal value. During this re-run the other core was busy with solver runs; in a separate test such concurrent load slowed a solve by about 6\%. Two measurements bound the remaining timing noise. The balanced cells of the third experiment reuse the seeds of the first, and in the re-run these three pairs of identical cells differed in total CPU time by at most 2.4\%. In the solver comparison of Section~6.8 of the main text, where every method was run three times per instance, the runs on instances with a median time of at least 1~s differed from that median by at most 13.4\% for Algorithm~1 of the main text, 15.6\% for model~(13) of the main text and 21.8\% for model~(14) of the main text (3.8\%, 2.3\% and 2.8\% at the median); shorter runs vary more. Probe counts are the exact measure of work.

\paragraph{Solver comparison (Section~6.8 of the main text)}

Two instances were processed at a time, each in a process pinned to its own core. On these 90 instances the times of Algorithm~1 of the main text in Table~10 of the main text exceed those of the re-run behind Table~6 of the main text and Table~S1 by 12.7\% at the median (interquartile range 9.5--17.7\%). The three methods in Table~10 of the main text were measured under identical conditions, so the time ratios are not affected, but absolute times should be compared within a table.

\section{Additional experimental tables}\label{S:tables}
Table~S1 gives the per-size CPU time statistics behind Table~6 of the main text; Table~S2 gives the per-cell probe counts behind Section~6.5 of the main text.

\begin{table}[htbp]
\centering
\caption{CPU time statistics per size, over all time ranges of the size.}
\begin{tabular}{lrrrc}
\toprule
Size & $N$ & Mean (s) & SD (s) & 95\% CI\\
\midrule
$10\times10$ & 250 & 0.127 & 0.047 & $[0.121,0.133]$\\
$10\times15$ & 250 & 0.243 & 0.078 & $[0.233,0.253]$\\
$20\times20$ & 250 & 1.219 & 0.403 & $[1.169,1.269]$\\
$25\times20$ & 250 & 1.796 & 0.578 & $[1.724,1.868]$\\
$30\times30$ & 250 & 4.322 & 1.396 & $[4.148,4.496]$\\
$30\times35$ & 250 & 5.572 & 1.557 & $[5.378,5.766]$\\
$40\times40$ & 250 & 11.530 & 3.648 & $[11.076,11.985]$\\
$40\times45$ & 250 & 14.614 & 3.946 & $[14.123,15.106]$\\
$50\times50$ & 250 & 26.350 & 6.281 & $[25.568,27.132]$\\
\bottomrule
\end{tabular}
\end{table}

\begin{table}[htbp]
\centering
\caption{Median number of probes without and with the incumbent cap, and median ratio (10 instances per cell).}
\small
\begin{tabular}{llrrrrrr}
\toprule
 & & \multicolumn{3}{c}{$a_i\in[10,50]$} & \multicolumn{3}{c}{$a_i\in[1,5]$}\\
\cmidrule(lr){3-5}\cmidrule(lr){6-8}
Size & Times & No cap & Cap & Ratio & No cap & Cap & Ratio\\
\midrule
$8\times8$ & $[1,40]$ & 751 & 196 & 0.256 & 722 & 188 & 0.231\\
$8\times8$ & $[10,20]$ & 648 & 218 & 0.355 & 510 & 160 & 0.328\\
$12\times12$ & $[1,40]$ & 4,138 & 740 & 0.185 & 3,630 & 718 & 0.212\\
$12\times12$ & $[10,20]$ & 1,925 & 579 & 0.264 & 2,438 & 625 & 0.242\\
$16\times16$ & $[1,40]$ & 9,789 & 1,614 & 0.168 & 8,230 & 1,406 & 0.166\\
$16\times16$ & $[10,20]$ & 6,505 & 1,466 & 0.214 & 4,770 & 1,100 & 0.212\\
$20\times20$ & $[1,40]$ & 19,163 & 3,099 & 0.163 & 15,097 & 2,433 & 0.160\\
$20\times20$ & $[10,20]$ & 9,496 & 2,040 & 0.216 & 11,016 & 2,492 & 0.210\\
\bottomrule
\end{tabular}
\end{table}

\section{Data of the examples of Jain et al.}\label{S:examples}
Table~S3 is an optimal allocation of the problem with idle phases allowed for Example~1 (Section~5.3 of the main text); Table~S4 is the phase partition of Example~2 transcribed from the published table (Section~5.4 of the main text).

\begin{table}[htbp]
\centering
\caption{An optimal allocation of \PTP{} for Example~1: triad $(06,32,00)$, value 38, phase flows $(F_1,F_2,F_3)=(19,61,0)$. Entries follow the phase markup of Table~4 of the main text; a dash denotes zero.}
\begin{tabular}{lccccccccc}
\toprule
 & $D_1$ & $D_2$ & $D_3$ & $D_4$ & $D_5$ & $D_6$ & $D_7$ & $D_8$ & $a_i$\\
\midrule
$S_1$ & \underline{6} & -- & -- & -- & \textbf{9} & -- & -- & -- & 15 \\
$S_2$ & -- & -- & -- & -- & -- & -- & \textbf{3} & \underline{7} & 10 \\
$S_3$ & \textbf{6} & -- & \textbf{6} & -- & -- & -- & -- & -- & 12 \\
$S_4$ & -- & -- & -- & \underline{6} & -- & -- & -- & -- & 6 \\
$S_5$ & -- & \textbf{5} & \textbf{7} & -- & -- & -- & \textbf{5} & -- & 17 \\
$S_6$ & -- & -- & -- & -- & -- & \textbf{20} & -- & -- & 20 \\
\midrule
$b_j$ & 12 & 5 & 13 & 6 & 9 & 20 & 8 & 7 & 80 \\
\bottomrule
\end{tabular}
\end{table}

\begin{table}[htbp]
\centering
\caption{Example~2 of~\cite{Jain2021} with the phase partition transcribed from the published table. Phase-1 times are underlined, Phase-2 times are in bold, and Phase-3 times are plain.}
\small
\begin{tabular}{lccccccccccc}
\toprule
 & $D_1$ & $D_2$ & $D_3$ & $D_4$ & $D_5$ & $D_6$ & $D_7$ & $D_8$ & $D_9$ & $D_{10}$ & $a_i$\\
\midrule
$S_{1}$ & \textbf{21} & \underline{10} & \underline{48} & 19 & \textbf{31} & \underline{58} & 43 & \textbf{39} & 60 & \underline{29} & 5 \\
$S_{2}$ & \underline{29} & 43 & \textbf{21} & 09 & \underline{35} & \textbf{12} & \underline{58} & 28 & \underline{10} & \textbf{45} & 12 \\
$S_{3}$ & 60 & \textbf{45} & \underline{17} & \textbf{12} & \underline{58} & 19 & 43 & \underline{17} & \textbf{45} & \underline{10} & 20 \\
$S_{4}$ & \underline{35} & \textbf{12} & 60 & \underline{29} & 28 & \textbf{31} & \underline{10} & 09 & \textbf{39} & \underline{48} & 14 \\
$S_{5}$ & \textbf{45} & \underline{35} & 09 & \textbf{61} & 19 & \underline{29} & \textbf{31} & \underline{58} & \underline{17} & \textbf{21} & 18 \\
$S_{6}$ & \underline{10} & \textbf{39} & \underline{29} & \textbf{39} & 43 & \textbf{61} & 28 & 60 & \textbf{12} & 19 & 21 \\
$S_{7}$ & 09 & \underline{48} & \textbf{31} & \textbf{12} & \underline{10} & 28 & \textbf{45} & \underline{58} & 43 & \underline{35} & 7 \\
$S_{8}$ & \textbf{12} & 28 & \underline{58} & \underline{35} & \textbf{21} & 09 & 60 & 19 & \underline{29} & \textbf{61} & 25 \\
$S_{9}$ & 19 & \underline{48} & 28 & \textbf{21} & 43 & \underline{10} & 09 & \textbf{45} & 43 & \underline{17} & 15 \\
$S_{10}$ & \textbf{31} & \underline{17} & \underline{35} & \underline{48} & \textbf{12} & 60 & \underline{29} & \textbf{21} & 09 & \textbf{39} & 6 \\
\midrule
$b_j$ & 18 & 7 & 14 & 6 & 21 & 12 & 20 & 12 & 6 & 27 & 143 \\
\bottomrule
\end{tabular}
\end{table}

\section{Filtering study: sampler, seeds and replicates}\label{S:filter}
This section gives the complete description of the filtering study of Section~6.9 of the main text, including the sampler, the seeds and the replicate samples.

For 2,000 random instances with $m,n\in\{3,4\}$, availabilities in $\{1,\dots,4\}$, times in $\{1,\dots,12\}$, phases sampled by shuffling the $mn$ links and cutting the sequence after $k_1$ and $k_1+k_2$ links, with $k_1$ uniform on $\{1,\dots,mn-2\}$ and $k_2$ uniform on $\{1,\dots,mn-k_1-1\}$ (so that Phase-1 receives about half of the links on average), and $\SNI\ne\emptyset$, we enumerated $\Shat$ and computed, for each phase $k$, the set $S^*_k$ of all triads that are proficient with respect to Phase-$k$, ties included. Every list that the algorithm of~\cite{Jain2021} can record when generating proficient triads with respect to Phase-$k$ is a subset of $S^*_k$, whatever its tie-breaking, by Theorem~4.5 of~\cite{Jain2021}, whose proof is the argument in the proof of Proposition~6 of the main text with $\Shat$ in place of $\SNI$; so if discarding the idle triads of $S^*_k$ does not give $\OPTNI$, filtering fails for every implementation, and the frequency of this event is a lower bound on the failure rate. We also ran the algorithm of~\cite{Jain2021} itself, with exact auxiliary solutions and random tie-breaking. Filtering fails, whatever the implementation, with respect to at least one phase in 1,597 instances (79.8\%) and with respect to all three in 525 (26.2\%); for the algorithm as run, the numbers are 1,671 (83.5\%) and 965 (48.2\%) (Table~11 of the main text). Among the 1,707 instances in which non-idleness binds, the lower bound for Phase-1 is 64.7\%. The rates of Table~11 of the main text are the frequencies observed in this sample, which was generated with seeds 2027 (its first 1,400 instances) and 2026 (the remaining 600) and is reproduced exactly from them; the intervals of the table express the sampling error of these frequencies. Twelve further samples of 2,000 instances from the same distribution gave Phase-1 lower bounds between 56.4\% and 61.4\%, all above the 55.2\% of this sample.\footnote{The twelve samples use seeds 100--111, and the eight samples with independent phase labels seeds 5001--5008. For them only the sets $S^*_k$ were computed, so their instance streams differ from that of the sample of Table~11 of the main text.} Pooled, they estimate the Phase-1 lower bound under this distribution at 59.2\% (14,201 of 24,000 instances; 95\% Wilson interval 58.5--59.8\%) and the lower bound for at least one phase at 81.9\% (81.4--82.4\%). The Phase-1 interval of Table~11 of the main text does not contain this estimate, so the table describes this sample, and its Phase-1 rate is lower than the estimated rate under the distribution. The rates also depend on the sampling scheme: when each link instead receives an independent uniform phase label, eight samples of 2,000 instances gave Phase-1 lower bounds between 40.4\% and 42.9\% (41.8\% pooled), and lower bounds for at least one phase between 67.2\% and 69.8\%.

The sampler shuffles the $mn$ links and cuts the sequence after $k_1$ and $k_1+k_2$ links, with $k_1$ uniform on $\{1,\dots,mn-2\}$ and $k_2$ uniform on $\{1,\dots,mn-k_1-1\}$; Phase-1 thus receives $(mn-1)/2$ links in expectation, while its share is uniform over $1,\dots,mn-2$. The archive reproduces the sample of Table~11 of the main text byte for byte from the seeds 2027 and 2026, and recomputes all twelve replicates and the eight samples with independent phase labels.

\section{HiGHS observations and reproduction}\label{S:highs}
This section documents the observations summarized in the footnote to Section~6.7 of the main text.

\paragraph{Summary} Before the timing runs, both formulations were also checked against Algorithm~1 of the main text on 3,300 further random instances with $m,n\le8$, five of which have no non-idle solution. In that check, with the settings of Section~6.8 of the main text, HiGHS 1.15.1 twice reported as optimal a value that was not optimal, once for each formulation: 24 where the optimum is 19, and 14 where it is 11. Both optima were established independently of Algorithm~1 of the main text and of the solver by an exact enumeration over phase thresholds, and for each saved model a feasible solution with the optimal value was checked against all its constraints. With presolve disabled, both runs returned the optimum. When the two saved models were solved again on another computer, the second discrepancy (14 instead of 11) recurred with HiGHS 1.13.0 to 1.15.1, but not with 1.11.0 or 1.12.0. The first (24 instead of 19) did not recur with the default options. With the option \texttt{random\_seed=15}, however, HiGHS 1.13.0 to 1.15.1 reported a different non-optimal value, 28, as optimal for the same saved model, whose optimum is 19; with presolve disabled they returned 19. We have reported this to the HiGHS developers as issue~\#3433 (\url{https://github.com/ERGO-Code/HiGHS/issues/3433})~\cite{HiGHSissue3433}. Agreement with the solver is therefore supporting evidence; the comparison with complete enumeration is the primary check.

\paragraph{Historical pre-check} The two discrepancies described in the summary arose in a cross-check of 3,300 random instances (300 with $m,n\le6$ and 3,000 with $m,n\le8$) against the algorithm of the main text, with HiGHS 1.15.1 on one thread and relative gap zero. The logged records are those of instance $r=292$ (direct model, $5\times6$): Optimal 24, optimum 19; and of instance $r=1404$ (threshold model, $8\times2$): Optimal 14, optimum 11. The saved threshold model omits the constant $\sum_kd_{k,1}=5$ of its objective, so its values are $9$ and $6$. With presolve disabled both runs returned the optimum. The platform of these runs was not recorded.

\paragraph{Reproduction from the saved model files} Both saved files were solved again on macOS~26.6 (arm64) with Python~3.13.5, each HiGHS release installed from PyPI in its own virtual environment, one thread, relative gap zero and presolve on (Table~\ref{tab:S5}); with presolve off every run returned the optimum. The 154 saved solution vectors were checked against the model constraints. The recorded decimal vectors satisfy feasibility within tolerance $10^{-6}$, and 129 are exactly feasible as recorded. After the documented reconstruction of integer and continuous values, all 154 are feasible in exact rational arithmetic. Thus the discrepancies are optimality claims for suboptimal feasible points.

\begin{table}[htbp]
\centering
\caption{Objective reported as Optimal (equal to the dual bound) for the saved models, presolve on, macOS arm64. Optima: 19 (direct model) and 6 (threshold model file).}\label{tab:S5}
\small
\begin{tabular}{@{}lccc@{}}
\toprule
HiGHS & \begin{tabular}[b]{@{}c@{}}threshold model,\\ default seed\end{tabular} & \begin{tabular}[b]{@{}c@{}}direct model,\\ default seed\end{tabular} & \begin{tabular}[b]{@{}c@{}}direct model,\\ \texttt{random\_seed=15}\end{tabular}\\
\midrule
1.11.0, 1.12.0 & 6 & 19 & 19\\
\begin{tabular}[t]{@{}l@{}}1.13.0, 1.13.1, 1.14.0,\\ 1.15.0, 1.15.1\end{tabular} & \textbf{9} & 19 & \textbf{28}\\
\begin{tabular}[t]{@{}l@{}}development branch,\\ commit f144f50\end{tabular} & 6 & 19 & \textbf{28}\\
\bottomrule
\end{tabular}
\end{table}

\paragraph{Issue \#3433 and follow-up} The direct-model failure with \texttt{random\_seed=15} was reported as HiGHS issue~\#3433~\cite{HiGHSissue3433} on 16 October 2026. On 17 October 2026 a HiGHS developer reported that it did not reproduce on the latest development branch and closed the issue. Source builds of the commits f144f50 and 8fd6a9a (the latest commit on 17 October 2026) and of release 1.15.1 were then tested on arm64 and, under Rosetta~2, on x86-64 (Table~\ref{tab:S6}). Seed~15 gives the optimum on the x86-64 builds of both development commits, but the seeds 124 and 267 give 28 on both architectures, and seeds 11 and 127 give 24 on x86-64. The failure therefore persisted on the latest development branch at that date; which seeds trigger it depends on the platform and the version. Linux was not tested.

\begin{table}[htbp]
\centering
\caption{Objective reported as Optimal (equal to the dual bound) for the saved direct model (optimum 19), presolve on, for selected values of \texttt{random\_seed}. With presolve off every run returned 19.}\label{tab:S6}
\small
\begin{tabular}{llcccccc}
\toprule
Build & Architecture & 0 & 11 & 15 & 124 & 127 & 267\\
\midrule
release 1.15.1 (PyPI) & arm64 & 19 & 19 & \textbf{28} & \textbf{28} & 19 & 19\\
commit f144f50 & arm64 & 19 & 19 & \textbf{28} & \textbf{28} & 19 & \textbf{28}\\
commit 8fd6a9a & arm64 & 19 & 19 & \textbf{28} & \textbf{28} & 19 & \textbf{28}\\
release 1.15.1 (PyPI) & x86-64 & 19 & 19 & \textbf{28} & \textbf{28} & \textbf{23} & 19\\
release 1.15.1 (source) & x86-64 & 19 & 19 & \textbf{28} & \textbf{28} & \textbf{23} & 19\\
commit f144f50 & x86-64 & 19 & \textbf{24} & 19 & \textbf{28} & \textbf{24} & \textbf{28}\\
commit 8fd6a9a & x86-64 & 19 & \textbf{24} & 19 & \textbf{28} & \textbf{24} & \textbf{28}\\
\bottomrule
\end{tabular}
\end{table}

The model files, feasible solutions, scripts, solver logs and environment records are in the directory \texttt{highs\_reproduction/} of the archive.

\FloatBarrier
\section{Two further instances for Remark~12 of Kaur, Dahiya and Verma}\label{S:remark12}
This section gives the two $2\times2$ instances referred to in Section~7.2 of the main text.

Both instances have $m=n=2$, $a=b=(2,2)$ and per-unit costs $c=\begin{pmatrix}1&2\\1&1\end{pmatrix}$. Every $X\in\Shat$ has the form $x_{11}=x_{22}=s$, $x_{12}=x_{21}=2-s$ with $s\in\{0,1,2\}$; denote it $X_s$. Its cost is $4$, $5$ and $6$ for $s=2,1,0$. Remark~12 of~\cite{Kaur2026arxiv} speaks of blocking Phase-II links in both of its cases. We read it as blocking links of the phase that carries all the flow, since blocking Phase-II links cannot make an empty Phase-II carry flow. Under the literal reading the cost-minimal solution remains $X_2$ in both instances below, so the re-solved problem never provides non-zero times of both phases.

\paragraph{Instance S6.1} Let $t=\begin{pmatrix}2&1\\1&2\end{pmatrix}$, with link $(1,2)$ in Phase-II and the other links in Phase-I. Then $X_2$ has $(T_1,T_2)=(2,0)$ and leaves Phase-II empty, $X_1$ has $(2,1)$, and $X_0$ has $(1,1)$; the non-dominated points are $(3,5)$ and $(2,6)$. The cost-minimal solution $X_2$ leaves Phase-II empty, so Remark~12 of~\cite{Kaur2026arxiv} re-solves the problem with some Phase-I links of time at least~2 blocked, that is, with a non-empty subset of $\{(1,1),(2,2)\}$ blocked. $X_1$ uses both links, so every choice excludes it, and the algorithm reports only $(2,6)$: it misses the non-dominated point of minimum cost.

\paragraph{Instance S6.2} Let $t_{ij}=1$ for all links, with Phase-I $=\{(1,1),(2,2)\}$ and Phase-II $=\{(1,2),(2,1)\}$. Then $X_2$ leaves Phase-II empty, $X_0$ leaves Phase-I empty, and $X_1$, with $(T_1,T_2)=(1,1)$ and cost~5, is the only non-idle solution. The cost-minimal solution $X_2$ leaves Phase-II empty, and every choice of blocked links removes $(1,1)$ or $(2,2)$, both of which $X_1$ uses. Hence $X_1$ is never reached, and the algorithm cannot return the only non-idle solution; under the readings we implemented, it reports that none exists.

\section{Computational checks of the proofs}\label{S:checks}
This section reports exhaustive computational checks of results of Section~4.5 of the main text.

As a check on the proofs, exhaustive search confirmed Lemma~5 of the main text on 5,641 random 3DM instances with $q\le4$, part~(a) of Theorem~5 of the main text on 3,146 reduced instances (2,128 with and 1,018 without a perfect matching) and directly on 793 instances that already had unique pairs, part~(b) on 3,943 reduced instances, and Proposition~4 of the main text against complete enumeration on 1,000 random instances for each of $p=2,3,4$; in every case the number of probes stayed within (12) of the main text.

\section{Contents of the archive}\label{S:archive}
The archive \path{three_phase_ni_supplement_matcom.zip} (which unpacks into the directory \path{three_phase_ni_supplement/}) contains a README with exact commands and checksums for every file, and the following parts.
\begin{itemize}
\item \texttt{filtering/}: the original filtering experiment, regenerated byte for byte from its seeds, with a script that recomputes every number of Section~S4 and of the corresponding table of the main text.
\item \texttt{sensitivity/}: the twelve replicate samples and the eight samples with independent phase labels.
\item \texttt{correctness/}: the correctness check on 5,000 instances, with an independent implementation of the algorithm and two exhaustive enumerations.
\item \texttt{highs\_reproduction/}: the material of Section~S5.
\item \texttt{jain/}: an independent count of the auxiliary problems of the algorithm of Jain et al.\ on the six-problem example and on random instances.
\item \texttt{dispatch/}: the synthetic dispatch scenario of the main text, its verification by enumeration and by an implementation of the algorithm.
\item \texttt{original\_studies/}: the original source code and recorded outputs of the computational study of Section~6 of the main text, including the solver comparison of Section~6.8 of the main text, of the study of Remark~12 of~\cite{Kaur2026arxiv} in Section~7.2 of the main text and of the checks of Section~S7, with a guide that maps the files to the reported results. These are historical records; they were not rerun for the archive.
\item \texttt{proposition6\_reproduction/}: The aggregate log of the October 2026 reproduction of the computational check of Proposition~6, together with a README documenting the execution commands and provenance limitations. The original run records and individual instance-level results are unavailable.
\end{itemize}

\bibliographystyle{elsarticle-num}
\bibliography{references}